\RequirePackage{fix-cm}
\documentclass[smallextended]{svjour3}       
\smartqed  

\usepackage{marvosym}

\usepackage{graphicx}
\begin{document}

\title{Cram\'{e}r moderate deviation expansion for martingales with one-sided Sakhanenko's condition and its applications 
}

\titlerunning{Cram\'{e}r moderate deviations for martingales}        

\author{Xiequan Fan   \and
 Ion Grama  \and
   Quansheng~Liu
}

\authorrunning{X. Fan,  I. Grama  and  Q. Liu } 

\institute{X. Fan     \at
Center for Applied Mathematics,
Tianjin University, Tianjin 300072,  China\\
             \email{fanxiequan@hotmail.com}
           \and
           I. Grama \and Q. Liu \\
            Univ. Bretagne-Sud,  UMR 6205, LMBA, 56000 Vannes, France\\
              \email{ion.grama@univ-ubs.fr, quansheng.liu@univ-ubs.fr}
}
\date{Received: date / Accepted: date}

\maketitle

\begin{abstract}
We give  a Cram\'{e}r moderate deviation expansion for  martingales  with differences   having finite conditional moments
of order $2+\rho,  \rho \in (0,1],$ and finite  one-sided conditional exponential moments.
 The upper bound of the range of validity and the remainder of our expansion are both optimal.
Consequently, it leads to a ``half-side''  moderate deviation principle for martingales.
Moreover, applications  to quantile coupling inequality, $\beta$-mixing and $\psi$-mixing sequences are discussed.


\end{abstract}

 \section{Introduction}
 \setcounter{equation}{0}
Let $(\eta_{i})_{i\geq 1}$ be a  sequence of independent and identically distributed (i.i.d.) centered real random variables (r.v.s)
satisfying Cram\'{e}r's condition $\mathbf{E} \exp\{ c_{0}|\eta_{1}|\} <\infty,$ for some constant $c_{0}>0.$
Without loss of generality,   assume that $ \mathbf{E}\eta_{1}^2=1.$
Cram\'{e}r \cite{Cramer38} established an asymptotic expansion of the probabilities of moderate deviations for the partial sums $\sum_{i=1}^{n}\eta_{i}$, based on the powerful technique of conjugate distributions (see also Esscher \cite{Esscher32}).
The result  of Cram\'{e}r implies that  uniformly in $0\leq x =  o(n^{1/2}), $
\begin{equation}
\log \frac {\mathbf{P}(\sum_{i=1}^{n}\eta_{i}> x \sqrt{n})} {1-\Phi(x)}  =  O \bigg( \frac{1+x^3}{\sqrt{n}}\bigg) \ \ \mbox{as} \ \ n \rightarrow \infty,
\label{cramer001}
\end{equation}
where $\Phi(x)=\frac{1}{\sqrt{2\pi}}\int_{-\infty}^{x}\exp\{-t^2/2\}dt$ is the standard normal distribution function.
Cram\'{e}r type moderate deviations for sums of independent r.v.s have been obtained by many authors. See, for instance,  Feller \cite{Fl43}, Petrov \cite{Pe54}, Sakhanenko \cite{S85} and \cite{FGL17}.
We refer to the monographs of Petrov \cite{Petrov75}, Saulis and Statulevi\v{c}ius \cite{SS78}  and the references therein.

In this paper we are concerned with Cram\'{e}r moderate deviations for martingales.
When the martingale differences are bounded, we refer to Bose \cite{Bose86a,Bose86b}, Ra\v{c}kauskas \cite{Rackauskas90,Rackauskas95,Rackauskas97},   Grama and Haeusler \cite{GH00}.
Let $(\eta _i,\mathcal{F}_i)_{i=0,...,n}$ be a sequence of  square
integrable martingale differences defined on a probability space $(\Omega, \mathcal{F}, \mathbf{P})$, where $\eta _0=0$ and $\{\emptyset, \Omega\}=\mathcal{F}_0\subseteq ...\subseteq
\mathcal{F}_n\subseteq \mathcal{F}$.
Assume that there exist absolute constants  $H>0$ and $N\geq 0$ such that $\max_i\left| \eta_i\right| \leq H$ and
$\left| \sum_{i=1}^n \mathbf{E}[\eta _i^2|\mathcal{F}_{i-1}]-n \right| \leq N^2.$
Here and hereafter, the equalities and inequalities between random variables are understood in the $\mathbf P$-almost sure sense.
From the results in Grama and  Haeusler \cite{GH00}, it follows that
\begin{equation}
\log \frac {\mathbf{P}(\sum_{i=1}^{n}\eta_{i}> x \sqrt{n})} {1-\Phi(x)}  =  O \bigg( \frac{x^3}{\sqrt{n}}\bigg),
\label{cramer002}
\end{equation}
for all $\sqrt{\log n}\leq x =  o(n^{1/4}), n \rightarrow \infty,$
and that
\begin{equation}
\frac {\mathbf{P}\left( \sum_{i=1}^{n}\eta_{i} >x\sqrt{n}\right)} {1-\Phi \left( x\right)} =
1+o\big(1\big)
\label{cramer003}
\end{equation}
uniformly for   $0 \leq x = o\left( n^{1/6}\right), n \rightarrow \infty.$
In  \cite{FGL13}
the expansions (\ref{cramer002}) and (\ref{cramer003}) have been extended to the case of martingale differences  satisfying the conditional Bernstein condition:
\begin{equation}
\Big|\mathbf{E}[\eta_{i}^{k}  | \mathcal{F}_{i-1}] \Big| \leq \frac{1}{2}k!H^{k-2} \mathbf{E}[\eta_{i}^2 | \mathcal{F}_{i-1}] \ \ \ \ \mbox{for}\ \ k\geq 3\ \ \mbox{and} \ \ 1\leq i\leq n,
\label{Bernst cond}
\end{equation}
where $H$ is a positive absolute constant.
We note that the conditional Bernstein condition implies that the martingale differences  have finite two-sided  conditional exponential moments.

In this paper we extend the expansions (\ref{cramer002}) and (\ref{cramer003}) to the case of martingales with differences having finite $(2+\rho)$th moments, $\rho \in (0, 1],$ and finite  one-sided conditional exponential moments. Assume that there exist  constants  $L,  M>0$ and $N\geq 0$ such that
\begin{equation}
\mathbf{E} [|\eta_{i}|^{2+\rho}e^{L \,\eta_{i}^+ } | \mathcal{F}_{i-1} ] \leq M^{ \rho} \, \mathbf{E} [ \eta_{i} ^{2} | \mathcal{F}_{i-1} ]   \ \ \ \textrm{for all} \ 1\leq i\leq n  \label{s08fd}
\end{equation}
and
\begin{equation}
 \Big| \sum_{i=1}^n \mathbf{E}[\eta _i^2|\mathcal{F}_{i-1}]-n \Big| \leq N^2.
\end{equation}
It is easy to see that  the conditional  Bernstein  condition implies   (\ref{s08fd}) with $\rho=1$,
while  condition (\ref{s08fd}) generally does not imply the conditional  Bernstein  condition; see (\ref{fhlm56}) for an example.
In Theorem \ref{co0} of the paper, we prove that if $\rho \in (0, 1),$ then for all  $ 0 \leq x =o\left( n^{1/2}\right),$
\begin{equation}
\log \frac { \mathbf{P}\left( \sum_{i=1}^{n}\eta_{i}>x\sqrt{n}  \right) } { 1-\Phi \left( x\right)} =
 O\bigg(\frac {1+x^{2+\rho}}{n^{\rho/2}}\bigg)  \ \ \mbox{as} \ \ n \rightarrow \infty.
\label{cramer004}
\end{equation}
The expansion (\ref{cramer004}) can be regard as an extension of  (\ref{cramer002}). We would like to point out that
the range of  validity of  (\ref{cramer002}) has been enlarged  to  the classical  Cram\'{e}r's one, and therefore  is optimal. Moreover, it is worth mentioning that  (\ref{cramer004}) is new   even for  independent r.v.s.
The last expansion implies that  (\ref{cramer003}) holds
uniformly in the range  $0 \leq x = o\left( n^{\rho/(4+2\rho)}\right).$ 
We also show that  when $\rho =1,$ equality (\ref{cramer004}) holds for all $\sqrt{\log n} \leq x =  o(n^{1/2}),$  see Remark   \ref{rema01} for details.


%
%

The paper is organized as follows. Our main results for martingales are stated and discussed in Section \ref{sec2}. Applications  to quantile coupling inequality, $\beta$-mixing and $\psi$-mixing sequences
are  discussed in Section \ref{sec3.1}.
Proofs of the theorems and their preliminary lemmas are deferred to Sections \ref{sec3}-\ref{secend}.
The proofs of Theorem \ref{co0} and Lemma \ref{LEMMA4} are refinements  of Fan et al.\ \cite{FGL13}.
The applications of our results  are new, and therefore are of independent interest.

Throughout the paper, $c$ and $c_\alpha,$ probably supplied with some indices,
denote respectively a generic positive  constant and a generic positive constant depending only on $\alpha.$
Denote by $\xi^+=\max\{\xi, 0 \}$ the positive part of $\xi$.
\section{Main results}\label{sec2}
\setcounter{equation}{0}
Let $n\geq1,$ and let $(\xi _i,\mathcal{F}_i)_{i=0,...,n}$ be   a sequence of martingale differences, defined on some
 probability space $(\Omega ,\mathcal{F},\mathbf{P})$,  where $\xi
_0=0 $,  $\{\emptyset, \Omega\}=\mathcal{F}_0\subseteq ...\subseteq \mathcal{F}_n\subseteq
\mathcal{F}$ are increasing $\sigma$-fields and $(\xi _i)_{i=1,...,n}$ are allowed to depend on $n$. Set
\begin{equation}
X_{0}=0,\ \ \ \ \ X_k=\sum_{i=1}^k\xi _i,\quad \ \   k=1,...,n.  \label{xk}
\end{equation}
Let $\left\langle X\right\rangle $ be the conditional variance of the
martingale $X=(X_k,\mathcal{F}_k)_{k=0,...,n}:$%
\begin{equation}\label{quad}
\left\langle X\right\rangle _0=0,\ \ \ \ \ \left\langle X\right\rangle _k=\sum_{i=1}^k\mathbf{E}[\xi _i^2|\mathcal{F}
_{i-1}],\quad \ \   k=1,...,n.
\end{equation}
In the sequel we shall use the following conditions:

\begin{description}
\item[(A1)]  There exist a constant  $\rho \in (0, 1]$ and  positive numbers  $\varepsilon_n \in (0, \frac12]$  such that
\[
\mathbf{E} [|\xi_{i}|^{2+\rho}e^{\varepsilon_n^{-1} \xi_{i}^+ } | \mathcal{F}_{i-1} ] \leq \varepsilon_n^{ \rho} \, \mathbf{E} [ \xi_{i} ^{2} | \mathcal{F}_{i-1} ]   \ \ \ \textrm{for all} \ 1\leq i\leq n.
\]
\item[(A2)]  There exist   non-negative numbers  $ \delta_n\in [0, \frac12]$ such that
$ \left| \left\langle X\right\rangle _n-1\right| \leq  \delta_n^2$\ \, a.s.

\end{description}

Condition (A1) can be seen as a one-sided version of  Sakhanenko's condition \cite{S85}.
In the case of normalized sums of i.i.d.\ random variables,  conditions (A1) and (A2) are satisfied with $\varepsilon_n = O(\frac {1} {  \sqrt n})$ and $\delta_n = 0.$
In the case of martingales, $\varepsilon_n$ and $\delta_n$ usually are satisfying
$\varepsilon_n, \delta_n \rightarrow 0$ as $n\rightarrow \infty$.

Notice that condition (A1) implies that $\mathbf{E} [ e^{\varepsilon_n^{-1} \xi_{i}^+ } | \mathcal{F}_{i-1} ]$ must be finite, which means that the positive
part of the conditional distribution of $\xi_{i}/\varepsilon_n$ has an exponential moment, and therefore has conditional moments of any order. However, such an assumption
is not required for the negative part of the conditional distribution.  For the negative part of $\xi_{i}$, we assume a finite conditional moment
of order $2+\rho$. Thus,   condition  (A1)   does not imply the conditional Cram\'{e}r condition, because  $\mathbf{E}[e^{\varepsilon_n^{-1} |\xi_{i}| }| \mathcal{F}_{i-1}] $ may not exist.

Let us remark that if $\xi_{i} $ is bounded, say $ |\xi_{i}|  \leq \gamma_n,$  then   condition  (A1)
is satisfied with $\varepsilon_n=e^{1/\rho} \,\gamma_n$.
On the other hand, if $\xi_{i} $ satisfies
\begin{equation}\label{fhlm56}
 \xi_{i}   \leq \gamma_n\ \ \
\textrm{and}\ \ \ \mathbf{E} [|\xi_{i}|^{2+\rho}  | \mathcal{F}_{i-1} ] \leq \tau_n^{ \rho} \, \mathbf{E} [ \xi_{i} ^{2} | \mathcal{F}_{i-1} ] \ \ \ \textrm{for all}\    1\leq i\leq n,
\end{equation}
then   condition  (A1) is also satisfied with
$\varepsilon_n = \max \{\gamma_n,  e^{1/\rho}\tau_n\}$. Here we assume that $ 0<\gamma_n,  \tau_n \leq \frac12 e^{-1/\rho} .$


The following theorem gives a Cram\'{e}r  moderate deviation expansion  for martingales.
\begin{theorem}\label{co0}
Assume conditions (A1) and (A2).
\begin{itemize}
   \item[\emph{[i]}] If $\rho \in (0, 1)$, then there is a constant $\alpha >0$, such that
for all $0\leq x \leq \alpha  \varepsilon_n^{-1},$
\begin{equation} \label{mt1ie1}
\bigg|\ln \frac{\mathbf{P}(X_n>x)}{1-\Phi \left( x\right)} \bigg| \leq   c_{\alpha,\rho} \,  \bigg(  x^{2+\rho} \varepsilon_n^\rho   + x^2 \delta_n^2 + (1+ x) \left( \varepsilon_n^\rho  +  \delta_n \right)\bigg)  .
\end{equation}
   \item[\emph{[ii]}] If $\rho =1$, then there is a constant $\alpha >0,$ such that
for all $0\leq x \leq \alpha  \varepsilon_n^{-1}$,
\begin{equation}\label{mt1ie2}
\bigg|\ln \frac{\mathbf{P}(X_n>x)}{1-\Phi \left( x\right)} \bigg| \leq  c_{\alpha}  \bigg(  x^{3} \varepsilon_n   + x^2 \delta_n^2 + (1+ x) \left( \varepsilon_n |\ln \varepsilon_n| +  \delta_n \right)\bigg)  .
\end{equation}
 \end{itemize}
\end{theorem}


The term $\varepsilon_n |\ln \varepsilon_n|$ in (\ref{mt1ie2}) cannot be replaced by $\varepsilon_n$ under the stated conditions. Indeed,
Bolthausen \cite{Bo82}  showed that there exists
a sequence of martingale  differences satisfying $|\xi_i|\leq 2/\sqrt{n}$ and $\langle X\rangle_n=1$ a.s., such that for all $n$ large enough,
\begin{equation} \label{ghklm02}
\sup_{x \in \mathbf{R}}\Big|\mathbf{P}(X_n \leq x)  -\Phi \left( x\right) \Big| \frac{\sqrt{n}}{\log n}   \ \geq \  c ,
\end{equation}
where $c$  is a  positive constant and does not  depend on $n$. See also \cite{F19} for general $\varepsilon_n$.  If $\varepsilon_n |\ln \varepsilon_n|$ in (\ref{mt1ie2}) could be improved to $\varepsilon_n,$ then we can deduce the following Berry-Esseen bound
\begin{equation}\label{ghklm01}
\sup_{x \in \mathbf{R}}|\mathbf{P}(X_n \leq x)  -   \Phi \left( x\right)   | \leq c\, (\varepsilon_n   +  \delta_n) ,
\end{equation}
which would violate   Bolthausen's result (\ref{ghklm02}). Thus $\varepsilon_n |\ln \varepsilon_n|$ in (\ref{mt1ie2}) cannot be improved to  $\varepsilon_n$  even for bounded martingale differences.

If the martingale differences are  bounded $|\xi_i|\leq \varepsilon_n$ and satisfy  condition (A2),  Grama and  Haeusler \cite{GH00} proved the asymptotic expansion  (\ref{mt1ie2}) for all $x \in [0,  \alpha   \min\{\varepsilon_n^{-1/2}, \delta_n^{-1}\}]$. Now Theorem \ref{co0} holds for a   larger range
$x \in [0,  \alpha \varepsilon_n  ^{-1}]$ and  a much more general class of martingales.

The following corollary states that under  conditions (A1) and (A2), the tail probabilities $\mathbf{P}( X_n>x)$ can be uniformly approximated by the tail probabilities of the standard normal  random variable, when $x$ is in a certain reduced range.
\begin{corollary}\label{co2}
Assume conditions (A1) and (A2).
\begin{itemize}
   \item[\emph{[i]}] If $\rho \in (0, 1)$, then
for all $0\leq x =o(\min\{\varepsilon_n^{-\rho/(2+\rho)} , \delta_n^{-1}\} )$,
\begin{equation}
\bigg|\frac{\mathbf{P}(X_n>x)}{1-\Phi \left( x\right)}- 1  \bigg| \leq    c_{\rho}\,   \bigg(  x^{2+\rho} \varepsilon_n^\rho   +  (1+ x) \left( \varepsilon_n^\rho  +  \delta_n \right)\bigg).
\end{equation}

   \item[\emph{[ii]}] If $\rho =1$, then
for all $0\leq x =o(\min\{\varepsilon_n^{-1/3} , \delta_n^{-1}\} )$,
\begin{equation}
\bigg|\frac{\mathbf{P}(X_n>x)}{1-\Phi \left( x\right)}-1  \bigg| \leq    c \,   \bigg(  x^{3} \varepsilon_n    +  (1+ x) \left( \varepsilon_n|\ln \varepsilon_n|  +  \delta_n \right)\bigg).
\end{equation}
 \end{itemize}
In particular, this implies that
  \[
  \frac{\mathbf{P}(X_n>x)}{1-\Phi \left( x\right)}=1+ o(1)
  \]
holds uniformly for $0\leq x =o(\min\{\varepsilon_n^{-\rho/(2+\rho)} , \delta_n^{-1}\} )$ as $\max\{\varepsilon_n, \delta_n \}\rightarrow 0.$

\end{corollary}

The inequalities (\ref{mt1ie1}) and (\ref{mt1ie2}) together implies that  there is a constant $\alpha >0$ such that
for $\rho\in (0 ,1 ]$ and all $0\leq x \leq \alpha  \varepsilon_n^{-1},$
\begin{equation}\label{fsvdf1}
\bigg| \log \frac{\mathbf{P}(X_n>x)}{1-\Phi \left( x\right)} \bigg|  \leq   c_{\alpha} \bigg(  x^{2+\rho} \varepsilon_n^\rho   + x^2 \delta_n^2 + (1+ x) \left( \varepsilon_n^\rho |\ln \varepsilon_n| +  \delta_n \right)\bigg).
\end{equation}
By (\ref{fsvdf1}), we obtain
the following   moderate deviation principle (MDP) result.
\begin{corollary}\label{mdp}
Assume that conditions (A1) and (A2) are satisfied with $ \max\{\delta_n, \varepsilon_n\} \rightarrow 0$   as $n\rightarrow \infty$.
Let $a_n$ be any sequence of real numbers satisfying $a_n \rightarrow \infty$ and $a_n \varepsilon_n \rightarrow 0$
as $n\rightarrow \infty$. Then for each Borel set $B\subset [0, \infty)$,
\begin{eqnarray}
- \inf_{x \in B^o}\frac{x^2}{2} &\leq & \liminf_{n\rightarrow \infty}\frac{1}{a_n^2}\log \mathbf{P}\left(\frac{1}{a_n} X_n \in B \right) \nonumber\\
 &\leq& \limsup_{n\rightarrow \infty}\frac{1}{a_n^2}\log \mathbf{P}\left(\frac{1}{a_n} X_n \in B \right) \leq  - \inf_{x \in \overline{B}}\frac{x^2}{2} \, , \label{ffhmts}
\end{eqnarray}
where $B^o$ and $\overline{B}$ denote the interior and the closure of $B$ respectively.
\end{corollary}

Since (\ref{ffhmts}) may not hold for all Borel set $B\subset (-\infty,0],$  inequality (\ref{ffhmts}) does not imply the usual MDP, but it can be seen as a  ``half-side" MDP.

Similar MDP results for martingales can be found in Dembo \cite{D96}, Gao \cite{G96} and Djellout \cite{D02}. For the most recent work on MDP for martingales with the conditional Cram\'{e}r condition and  the assumption that $  \mathbf{E}[ \xi_i^2 | \mathcal{F}_{i-1}] = 1/n$ a.s.\ for all $i.$, we refer to Eichelsbacher and L\"{o}we \cite{PL17} where the authors established a MDP result via Lindeberg's method.

\begin{remark}  \label{rema01}
The sequence of martingale differences  $(\xi _i,\mathcal{F}_i)_{i=0,...,n} $ discussed so far is standardized. For a general sequence of martingale differences
$(\eta_i,\mathcal{F}_i)_{i\geq1}$, one can
restate the conditions (A1) and (A2) as below.
\begin{description}
\item[(A1$'$)] There exist three positive constants  $\rho \in (0, 1], K $ and $L$  such that
\[
\mathbf{E} [|\eta_{i} |^{2+\rho}e^{ K  \eta_{i}^+ } | \mathcal{F}_{i-1} ] \leq L^{ \rho} \, \mathbf{E} [ \eta_{i}^{2} | \mathcal{F}_{i-1} ]   \ \ \ \textrm{for all} \ 1\leq i\leq n;
\]
\item[(A2$'$)] There exists a constant  $N\geq 0$   such that $$ \Big|\sum_{i=1}^{n}\mathbf{E}[\eta^2_{i} | \mathcal{F}_{i-1}] -n   \Big| \leq  N^2\ \ \ \  \textrm{a.s.}$$
\end{description}
Under conditions (A1$'$) and (A2\,$'$),  the  inequalities  (\ref{mt1ie1})-(\ref{ffhmts}) remain valid for
\begin{equation}
W_n = \sum_{i=1}^n \frac{\eta_i}{\sqrt{n}   }
\end{equation}
instead of $X_n,$ with $\varepsilon_n=n^{-1/2}\max\{K, L\}$ and $\delta_n=n^{-1/2}  L.$
\end{remark}

\section{Applications}\label{sec3.1}
\setcounter{equation}{0}
\subsection{Quantile coupling inequality}

Thanks to the work of Mason and Zhou \cite{MZ12}, it is  known that the Cram\'{e}r moderate deviation expansion
can be applied to establishing quantile coupling inequalities.   When the martingale differences  are  bounded,
a quantile coupling inequality  has been established by Mason and Zhou, see Corollary 2 of \cite{MZ12}.
Here, we give a generalization of the inequality of  Mason and Zhou \cite{MZ12}.

 Let $(W_n)_{n\geq 1}$ be a sequence of random variables and for each integer $n\geq1$, and let
$$F_n(x)=\mathbf{P}(W_n \leq x), \ \ x \in \mathbf{R},$$
denote the cumulative distribution function of $W_n.$  Its \textit{quantile function} is defined
by
$$H_n(s)=\inf\{ x : F_n(x)\geq s\}, \ \ s \in (0, 1).$$
Let $Z$ denote a standard normal random variable. Since $\Phi(Z)=_d U$ the uniformly distribution random variable,
then it is obvious that for each integer $n\geq1,$
$$ H_n( \Phi(Z))=_d W_n, $$ where $ =_d$ stands for equivalent in distribution.  For this reason,
we define
\begin{eqnarray}\label{fsffs}
W_n = H_n( \Phi(Z)) .
\end{eqnarray}
By Theorem \ref{co0}, we prove  the following quantile inequality.
\begin{theorem}\label{th3.0}
Let $(\eta _i,\mathcal{F}_i)_{i\geq1}$ be a sequence of martingale differences satisfying the following conditional Sakhanenko condition
\[
\mathbf{E} [|\eta_{i}|^{3}e^{K |\eta_{i}| } | \mathcal{F}_{i-1} ] \leq  L \, \mathbf{E} [ \eta_{i} ^{2} | \mathcal{F}_{i-1} ],\ \ \ \ \  i \geq 1,
\]
and
\[
\Big| \sum_{i=1}^n \mathbf{E} [ \eta_{i}^2 | \mathcal{F}_{i-1} ]  -n   \Big|\leq  M \ \  \textrm{a.s.},
\]
where  $\rho \in (0, 1],$ $K,$ $L$ and $M$ are positive constants. 
Assume that $W_n=_d\sum_{i=1}^n \eta_{i}/\sqrt{n}$ and $W_n$ is defined as in (\ref{fsffs}). There there exist constants $\alpha>0$ and $D>0$
and an integer $n_0$   such that
whenever $n\geq n_0$ and
\begin{eqnarray}\label{fg2sf}
 |W_n| \leq  \alpha \sqrt{ n},
\end{eqnarray}
we have
\begin{eqnarray}\label{fg235}
 \sqrt{n}|W_n-Z|/\ln n \leq 2 D (W_n^2 +1) \ \ \    \textrm{a.s.}
\end{eqnarray}
Furthermore, there exist two positive constants $C$ and $\lambda$ such that whenever $n\geq n_0$, we have for all $x\geq0,$
 \begin{eqnarray}\label{fg2sfs5}
 \mathbf{P}\Big( \sqrt{n}|W_n-Z|/\ln n > x \Big) \leq C \exp\Big\{ - \lambda \, x\, \Big\}.
\end{eqnarray}
\end{theorem}

When the martingale differences  are \textit{bounded}, Mason and Zhou \cite{MZ12} proved that (\ref{fg235}) holds whenever
$ |W_n| \leq  \alpha \sqrt[4]{ n}.$ Notice that the bounded martingale differences satisfy the conditional Sakhanenko condition.
Moreover, the range $ |W_n| \leq  \alpha \sqrt[4]{ n}$ has been extended to a much larger one  $|W_n| \leq  \alpha \sqrt { n}$  in our theorem.

\subsection{$\beta$-mixing sequences}\label{section02}

Let $(\eta_i )_{i\geq1}$ be a random process that may be non-stationary.
Write $S_{k,m}=\sum_{i=k+1}^{k+m}\eta_i.$
Assume that
there exists a constant  $\rho \in  (0, 1) $ such that
\begin{eqnarray}
 \mathbf{E}\eta_i =0   \ \ \ \textrm{for all} \   i , \label{cosf3.10}
\end{eqnarray}
\begin{eqnarray}
 \mathbf{E}|S_{k,m}|^{2+\rho} \leq m^{1+\rho/2}c_1^{2+\rho} , \label{cosf3.11}
\end{eqnarray}
and
\begin{eqnarray}\label{cosdf3.12}
 \mathbf{E}S_{k,m}^2 \geq c_2^2 m \ \ \ \  \textrm{for all} \  k\geq0, m\geq 1.
\end{eqnarray}
 Let $\mathcal{F}_{j} $ and $\mathcal{F}_{j+k}^{\infty}$
be $\sigma$-fields generated respectively by $(\eta_i)_{i \leq j}$ and $ (\eta_i)_{i \geq j+k}.$
 We say that $(\eta_i )_{i\geq1}$
is $\beta$-\textit{mixing}   if
$$\beta(n)=:  \sup_j \mathbf{E} \sup\{ \big|\mathbf{P}( B |\mathcal{F}_{j} ) -\mathbf{P}(B) \big|:\ B \in \mathcal{F}_{j+n}^{\infty}    \} \rightarrow \infty,\ \ \ \   \ \ n\rightarrow \infty.$$
Assume that there exist positive numbers $a_1, a_2$ and $\tau$ such that
\begin{eqnarray}\label{sgkm01}
\beta(n) \leq a_1 \exp\{-a_2 n^{\tau}\}.
\end{eqnarray}
By Theorem 4.1 of Shao and Yu \cite{SY96}, it is known that (\ref{cosf3.11})  is implied by the condition that $\mathbf{E}|\eta_{i}|^{2+\rho'} \leq  c_1^{2+\rho'}$ for a constant $\rho' >\rho$.

Set $ \alpha \in (0 , \frac12).$ Let $m=\lfloor n^\alpha \rfloor$ and $k=\lfloor n/(2m) \rfloor$ be respectively the integers part of $n^\alpha$ and $n/(2m)$.
Let
$$Y_j=\sum_{i=1 }^{ m  }\eta_{2m(j-1)+i} \ \ \ \ \ \  \textrm{and} \ \ \ \ \ \
  S_n=\sum_{j=1}^k Y_j.$$
Note that $S_n$ is an interlacing sum of $(\eta_i )_{i\geq1}$, and that $ \textrm{Var}  (S_n) = \mathbf{E}S_n^2.$

\begin{theorem} \label{tghbt}
Assume conditions   (\ref{cosf3.10})-(\ref{sgkm01}). Suppose that $ \eta_i \leq c_3$ for all $i.$
Then for all  $0  \leq x   =o( \min\{ n^{\frac12 - \alpha}, \,  n^{\alpha \, \tau /2} \}  ),$
\begin{equation}
\bigg | \ln \frac{\mathbf{P}(S_n/\sqrt{ \mathbf{E}S_n^2} >x)}{1-\Phi \left( x\right)} \bigg| \leq  \,  c_\rho \,   \frac{ (1+x)^{2+\rho }}{n^{\rho(\frac12 - \alpha)}  }    .
\end{equation}
In particular, we have
\begin{equation}
\frac{\mathbf{P}(S_n/\sqrt{ \mathbf{E}S_n^2}  >x)}{1-\Phi \left( x\right)}= 1+ o(1)
\end{equation}
uniformly for $0\leq x=o( \min\{ n^{ \rho(1-2 \alpha)/(4+2\rho)  }, \,  n^{\alpha \, \tau /2} \}  )$.
\end{theorem}

For a counterpart of Theorem  \ref{tghbt} for interlacing self-normalized   sums  $  W_n=S_n / \sqrt{\sum_{j=1}^kY_j^2},$ we refer to Chen et al.\ \cite{CSW16}.

The following MDP result is a consequence of the last theorem with $\alpha=1/(2+\tau)$.
\begin{corollary} Assume the conditions of Theorem  \ref{tghbt}.
Let $a_n$ be any sequence of real numbers satisfying $a_n \rightarrow \infty$ and $a_n  n^{-\tau/(2\tau+4)} \rightarrow 0$
as $n\rightarrow \infty$. Then for each Borel set $B \subset [0, \infty)$,
\begin{eqnarray*}
- \inf_{x \in B^o}\frac{x^2}{2} &\leq & \liminf_{n\rightarrow \infty}\frac{1}{a_n^2}\log \mathbf{P}\bigg(\frac{1}{a_n} \frac{S_n}{\sqrt{\mathbf{E}S_n^2}}   \in B \bigg) \nonumber\\
 &\leq& \limsup_{n\rightarrow \infty}\frac{1}{a_n^2} \log \mathbf{P}\bigg(\frac{1}{a_n} \frac{ S_n}{\sqrt{\mathbf{E}S_n^2}} \in B \bigg) \leq  - \inf_{x \in \overline{B}}\frac{x^2}{2} \,,
\end{eqnarray*}
where $B^o$ and $\overline{B}$ denote the interior and the closure of $B$ respectively.
\end{corollary}

\subsection{$\psi$-mixing sequences}
Recall the notations in Section \ref{section02}.
We say that $(\eta_i )_{i\geq1}$
is  $\psi$-\textit{mixing} if
\begin{eqnarray}\label{sgkm02}
 \psi(n)=:  \sup_j   \sup_B\{ \big|\mathbf{P}( B |\mathcal{F}_{j} ) -\mathbf{P}(B) \big|/\mathbf{P}(B):\ B \in \mathcal{F}_{j+n}^{\infty}    \} \rightarrow 0,\ \ \ \   \ \ n\rightarrow \infty.
\end{eqnarray}
Set $ \alpha \in (0, \frac12).$ Let $m=\lfloor n^\alpha \rfloor$ and $k=\lfloor n/(2m) \rfloor$ be respectively the integers part of $n^\alpha$ and $n/(2m),$
and let
$$Y_j=\sum_{i=1 }^{ m  }\eta_{2m(j-1)+i} \ \ \ \ \ \  \textrm{and} \ \ \ \ \ \
  S_n=\sum_{j=1}^k Y_j$$
 as in Section \ref{section02}. 

Denote
\begin{eqnarray}
  \tau_n^2=  \psi(m ) + n \psi^2(m )   +   k\psi^{1/2} (m )  .  \label{detan}
\end{eqnarray}
We have the following Cram\'{e}r  moderate deviations for $\psi$-mixing sequences.
\begin{theorem}\label{th3.3}
Assume conditions (\ref{cosf3.10})-(\ref{cosdf3.12}) with $\rho \in (0, 1]$. Suppose that $ \eta_i \leq c_3$ for all $i,$
and that $ \tau_n \rightarrow 0$ as $n\rightarrow \infty$.
\begin{itemize}
   \item[\emph{[i]}] If $\rho \in (0, 1)$, then
for all  $0\leq x  =o( n^{\frac1 2-\alpha} ),$
\begin{equation}\label{ggdfg01}
\bigg|\ln \frac{\mathbf{P}(S_n/\sqrt{ \mathbf{E}S_n^2} >x)}{1-\Phi \left( x\right)} \bigg| \leq   c_{\rho} \,  \bigg(  \frac{ x^{2+\rho }}{n^{\rho(\frac12 - \alpha)}  }  + x^2 \tau_n^2 + (1+ x) \Big( \frac{1}{n^{\rho(\frac12 - \alpha)}  } +  \tau_n  \Big)\bigg) .
\end{equation}

   \item[\emph{[ii]}] If $\rho =1$, then
for all  $0\leq x  =o( n^{\frac1 2-\alpha} ),$
\begin{equation}\label{ggdfg02}
\bigg|\ln \frac{\mathbf{P}(S_n/\sqrt{ \mathbf{E}S_n^2} >x)}{1-\Phi \left( x\right)} \bigg| \leq  c \, \bigg(\frac{  x^3 }{n^{ \frac12 - \alpha} }+ x^2\tau_n^2 + (1+ x) \Big( \frac{ |\ln n |}{n^{\frac12 - \alpha}}  +  \tau_n \Big) \bigg)  .
\end{equation}
 \end{itemize}
In particular, if
 \begin{eqnarray}\label{tphin}
\psi(n) = O\big( n^{- ( 2+\rho)(1-\alpha)/  \alpha  } \big) ,
\end{eqnarray}
then
\begin{equation}
\displaystyle \tau_n = O(n^{-\rho(\frac12 - \alpha)}  ) \ \ \ \ \ \textrm{and} \ \ \ \ \  \frac{\mathbf{P}(S_n/\sqrt{ \mathbf{E}S_n^2}  >x)}{1-\Phi \left( x\right)}= 1+ o(1)
\end{equation}
uniformly for $0\leq x=o( n^{ \rho(1-2 \alpha)/(4+2\rho)  }).$
\end{theorem}

In the independent case, we have $\psi(n)=0$ and $\tau_n=0$. Let $\alpha \rightarrow 0.$   Then (\ref{ggdfg01}) and (\ref{ggdfg02})
 recover the optimal  range of validity, that is $0\leq x =  o(n^{1/2}).$

The following MDP result is a consequence of the last theorem.
\begin{corollary} Assume the conditions of Theorem \ref{th3.3}.
Let $a_n$ be any sequence of real numbers satisfying $a_n \rightarrow \infty$ and $a_n/ n^{  \frac1  2  -\alpha } \rightarrow 0$
as $n\rightarrow \infty$. Then for each Borel set $B \subset [0, \infty)$,
\begin{eqnarray*}
- \inf_{x \in B^o}\frac{x^2}{2} &\leq & \liminf_{n\rightarrow \infty}\frac{1}{a_n^2}\log \mathbf{P}\bigg(\frac{1}{a_n} \frac{S_n}{\sqrt{\mathbf{E}S_n^2}}   \in B \bigg) \nonumber\\
 &\leq& \limsup_{n\rightarrow \infty}\frac{1}{a_n^2} \log \mathbf{P}\bigg(\frac{1}{a_n} \frac{ S_n}{\sqrt{\mathbf{E}S_n^2}} \in B \bigg) \leq  - \inf_{x \in \overline{B}}\frac{x^2}{2} \,,
\end{eqnarray*}
where $B^o$ and $\overline{B}$ denote the interior and the closure of $B$ respectively.
\end{corollary}


\section{Preliminary lemmas}\label{sec3}
\setcounter{equation}{0}
Assume condition (A1).   For any real  $\lambda \in [0, \ \varepsilon_n^{-1} ],$
define the  exponential multiplicative martingale  $Z(\lambda
)=(Z_k(\lambda ),\mathcal{F}_k)_{k=0,...,n},$ where
\[
Z_k(\lambda )=\prod_{i=1}^k\frac{e^{\lambda \xi _i}}{\mathbf{E}[e^{\lambda \xi _i}|
\mathcal{F}_{i-1}]},\quad k=1,...,n,\quad Z_0(\lambda )=1.  \label{C-1}
\]
Then for each $k=1,...,n,$ the random variable $Z_k(\lambda
) $ defines a probability density on $(\Omega ,\mathcal{F},\mathbf{P}).$ This allows us to introduce  the  conjugate probability measure  $\mathbf{P}_\lambda $ on $(\Omega ,%
\mathcal{F})$ defined by
\begin{equation}
d\mathbf{P}_\lambda =Z_n(\lambda )d\mathbf{P}.  \label{f21}
\end{equation}
Denote by $\mathbf{E}_{\lambda}$ the expectation with respect to $\mathbf{P}_{\lambda}$.
For all $i=1,\dots,n$, let
\[
\eta_i(\lambda)=\xi_i - b_i(\lambda)\ \ \ \ \ \ \ \ \  \textrm{where} \  \ \   b_i(\lambda)=\frac{\mathbf{E}[\xi
_ie^{\lambda \xi _i}|\mathcal{F}_{i-1}]}{\mathbf{E}[e^{\lambda \xi _i}|\mathcal{F}%
_{i-1}]}.  \]
We thus have the following decomposition:
\begin{equation}
X_k=Y_k(\lambda )+B_k(\lambda ),\quad\quad\quad k=1,...,n, \label{xb}
\end{equation}
where $Y(\lambda )=(Y_k(\lambda ),\mathcal{F}_k)_{k=1,...,n}$ is the
conjugate martingale defined as
\begin{equation}\label{f23}
Y_k(\lambda )=\sum_{i=1}^k\eta _i(\lambda ),\quad\quad\quad k=1,...,n,
\end{equation}
and
$B(\lambda )=(B_k(\lambda ),\mathcal{F}_k)_{k=1,...,n}$ is the
drift process defined as
\[
B_k(\lambda )=\sum_{i=1}^kb_i(\lambda ),\quad\quad\quad k=1,...,n.
\]

In the proofs of theorem,  we need  a two-sided bound for the drift process $B_n(\lambda ).$
To this end, we prove the following   lemma.
\begin{lemma}\label{lemma2s}
If there exists an $s> 2,$  such that
\begin{equation}
\mathbf{E} [|\xi_{i}|^{s}e^{\varepsilon_n^{-1}\xi_{i}^+  } | \mathcal{F}_{i-1} ] \leq \varepsilon_n^{s-2} \, \mathbf{E} [\xi_{i}^{2} | \mathcal{F}_{i-1} ],
\end{equation}
then
\begin{equation}\label{ineq28}
\mathbf{E} [\xi_{i}^{2} | \mathcal{F}_{i-1} ] \leq \varepsilon_n^2.
\end{equation}
In particular, condition (A1) implies    (\ref{ineq28}).
\end{lemma}
\noindent\textbf{Proof.}
By Jensen's inequality, it is easy to see that
\begin{eqnarray*}
(\mathbf{E} [\xi_{i}^{2} | \mathcal{F}_{i-1} ])^{s/2} \leq \mathbf{E} [|\xi_{i}|^{s}  | \mathcal{F}_{i-1} ]  \leq  \mathbf{E} [|\xi_{i}|^{s}e^{\varepsilon_n^{-1}\xi_{i}^+  } | \mathcal{F}_{i-1} ]
\leq      \varepsilon_n^{s-2} \, \mathbf{E} [\xi_{i}^{2} | \mathcal{F}_{i-1} ].
\end{eqnarray*}
Thus $$(\mathbf{E} [\xi_{i}^{2} | \mathcal{F}_{i-1} ])^{s/2 -1} \leq   \varepsilon_n^{s-2},$$
 which implies (\ref{ineq28}).  \hfill\qed

Using the last   lemma, we establish a two-sided bound for the drift process $B_n(\lambda ).$
\begin{lemma}
\label{LEMMA-1-1} Assume conditions (A1) and (A2). Then  for all $0\leq\lambda\leq   \varepsilon_n^{-1},$
\begin{equation}\label{f25sa}
| B_n(\lambda ) - \lambda |\leq
 \lambda \delta_n^{2}+ c  \,  \lambda^{1+\rho} \varepsilon_n^{ \rho} .
\end{equation}
\end{lemma}

\noindent\textbf{Proof.}
Jensen's inequality and $\mathbf{E}[\xi _i|\mathcal{F}_{i-1}]=0$ imply that $\mathbf{E}[e^{\lambda \xi _i}|\mathcal{F}
_{i-1}]\geq 1, \lambda\geq0.$ Notice that
\[
\mathbf{E}[\xi_{i} e^{\lambda\xi_{i}} |\mathcal{F}_{i-1}]=\mathbf{E}\big[\xi_{i}(e^{\lambda\xi_{i}}-1)|\mathcal{F}_{i-1} \big]\geq 0,\ \ \ \ \  0\leq\lambda\leq   \varepsilon_n^{-1}.
\]
Using Taylor's expansion for $e^x$, we get
\begin{eqnarray}
B_n(\lambda ) & \leq & \sum_{i=1}^{n}\mathbf{E}[\xi_{i} e^{\lambda \xi_{i}} | \mathcal{F}_{i-1}]\nonumber\\
& = & \lambda\langle X\rangle_{n}+\sum_{i=1}^{n}\mathbf{E}\big[\xi_{i}(e^{\lambda \xi_{i}}-1- \lambda \xi_{i})\, \big| \mathcal{F}_{i-1} \big].\nonumber
\end{eqnarray}
Recall $\rho \in (0, 1].$
When $x\leq-1,$  by  Taylor's expansion, it is easy to see that $ \big| x(e^x-1-x) \big|  \leq \big| x(e^x-1) \big|+ x^2 \leq   2 |x|^{2+\rho} .$
When $x \in (-1 ,  1),$ again by  Taylor's expansion, we get  $|x(e^x-1-x)| \leq \frac12 |x|^3 e^{ x^+} \leq |x|^{2+\rho}e^{x^+}.$
When $x \geq 1,$ we have $|x(e^x-1-x)| \leq xe^x \leq x^{2+\rho}e^{x}.$
Thus, it holds
\begin{eqnarray}\label{f2s9s}
 |x(e^x-1-x)| \leq   2 |x|^{2+\rho}e^{x^+}, \ \ \   \  \  x \in \mathbf{R}.
\end{eqnarray}
By   inequality (\ref{f2s9s}),
we obtain for all $0\leq\lambda\leq   \varepsilon_n^{-1},$
\begin{eqnarray}
B_n(\lambda )
  &\leq&     \lambda\langle X\rangle_{n}+ 2 \lambda^{1+\rho} \sum_{i=1}^{n}\mathbf{E}[|\xi_{i}|^{2+\rho}e^{\lambda \xi_i^+ }  | \mathcal{F}_{i-1} ] \nonumber \\
  &\leq& \lambda\langle X\rangle_{n}+  2\lambda^{1+\rho} \sum_{i=1}^{n}\mathbf{E}[|\xi_{i}|^{2+\rho}e^{\varepsilon_n^{-1} \xi_i^+}  | \mathcal{F}_{i-1} ].\label{f26}
\end{eqnarray}
Condition (A2) implies that $\langle X \rangle_{n}  \leq 2.$
Combining (\ref{f26}), conditions (A1) and (A2) together, we  get the upper bound of $B_n(\lambda )$:
\[
B_n(\lambda ) \leq \lambda\langle X\rangle_{n}+ 2 \lambda^{1+\rho}     \varepsilon_n^\rho   \langle X \rangle_{n} \leq
\lambda +  \lambda \delta_n^{2} + 4\lambda^{1+\rho}  \varepsilon_n^\rho.
\]
When $x\leq -1,$   by Taylor's expansion, it is easy to see that $ \big|e^x-1-x -\frac12 x^2 \big| \leq  \big|e^x-1-x\big|+  \frac12 x^2   \leq    |x|^{2+\rho} .$
When $x \in (-1 ,  1),$ again by  Taylor's expansion, we get  $  \big|e^x-1-x -\frac12 x^2 \big| \leq  \frac16 |x|^3 e^{ x^+} \leq |x|^{2+\rho}e^{x^+}.$
When $x \geq 1,$ we have $\big|e^x-1-x -\frac12 x^2 \big| \leq  \big|e^x-1-x\big|+ \frac12 x^2   \leq x^{2+\rho}e^{x}.$
Thus, it holds
\begin{eqnarray}\label{inds35}
 \Big|e^x-1-x -\frac12 x^2 \Big| \leq     |x|^{2+\rho}e^{x^+}, \ \ \   \rho \in (0, 1]\  \textrm{and}\  x \in \mathbf{R}.
\end{eqnarray}
Using inequality (\ref{inds35}), condition  (A1) and Lemma \ref{lemma2s}, we have for all $0\leq\lambda\leq   \varepsilon_n^{-1},$
\begin{eqnarray}
\mathbf{E}[e^{\lambda \xi_{i}} | \mathcal{F}_{i-1} ] & = & 1 +  \frac12 \lambda^2  \mathbf{E} [\xi_{i}^2 | \mathcal{F}_{i-1}  ]+ \mathbf{E}\big[e^{\lambda \xi_{i}}-1 - \lambda \xi_{i} -\frac12 \lambda^2 \xi_{i}^2\big| \mathcal{F}_{i-1} \big] \nonumber \\
& \leq & 1+  \frac12 \lambda^2\, \mathbf{E}\big[\xi_{i}^2 | \mathcal{F}_{i-1} \big]+   \lambda^{2+\rho} \mathbf{E}\big[   |\xi_{i}|^{2+\rho} e^{\varepsilon_n^{-1} \xi_i^+} \big| \mathcal{F}_{i-1} \big] \nonumber \\
& \leq & 1+  \big(\frac12 \lambda^2\,+   \lambda^{2+\rho}\varepsilon_n^\rho \big)\mathbf{E}\big[\xi_{i}^2 | \mathcal{F}_{i-1} \big]  \label{sbopfsf} \\
& \leq & 1+  2 (\lambda\varepsilon_n )^2 .  \label{sbopdf}
\end{eqnarray}
By inequality (\ref{f2s9s}) and the fact $\langle X \rangle_{n}  \leq 2$,  we deduce that  for all $0\leq\lambda\leq   \varepsilon_n^{-1},$
\begin{eqnarray*}
   \sum_{i=1}^{n}\mathbf{E}[\xi_{i} e^{\lambda \xi_{i}} | \mathcal{F}_{i-1}]
  & = &  \lambda\langle X\rangle_{n}+\sum_{i=1}^{n}\mathbf{E}\big[\xi_{i}(e^{\lambda \xi_{i}}-1- \lambda \xi_{i})\, \big| \mathcal{F}_{i-1} \big] \nonumber\\
  & \geq &   \lambda\langle X\rangle_{n} -  2 \lambda^{1+\rho}   \sum_{i=1}^{n}\mathbf{E}[|\xi_{i}|^{2+\rho} e^{\varepsilon_n^{-1} \xi_i^+} | \mathcal{F}_{i-1} ]  \\
   & \geq &   \lambda\langle X\rangle_{n} -  2 \lambda^{1+\rho}\varepsilon_n^\rho  \langle X\rangle_{n}  \\
  & \geq &  \lambda -\lambda \delta_n^2 -  4  \,  \lambda^{1+\rho}  \varepsilon_n^\rho.
\end{eqnarray*}
The last inequality together with   (\ref{sbopdf}) imply the lower bound of $B_{n}(\lambda)$: for all $0\leq\lambda\leq   \varepsilon_n^{-1},$
\begin{eqnarray*}
  B_n(\lambda )
  & \geq &  \Big(\lambda -\lambda \delta_n^2 -  4  \,  \lambda^{1+\rho}  \varepsilon_n^\rho \Big)\Big( 1 +  2\,(\lambda\varepsilon_n)^2 \Big)^{-1} \\
  & \geq &  \lambda - \lambda \delta_n^{2} - 6 \lambda^{1+\rho}  \varepsilon_n^\rho,
\end{eqnarray*}
where the last line follows from the following inequality
\begin{eqnarray*}
\lambda -\lambda \delta_n^2 -  4\,    \lambda^{1+\rho}  \varepsilon_n^\rho &\geq& \lambda -\lambda \delta_n^2   - ( 6 -     2(\lambda  \varepsilon_n)^{2-\rho}  )  \lambda^{1+\rho}  \varepsilon_n^\rho   \\
&\geq& \Big(\lambda - \lambda \delta_n^{2} - 6 \lambda^{1+\rho}  \varepsilon_n^\rho\Big)\Big(1 + 2 (\lambda\varepsilon_n)^2 \Big).
\end{eqnarray*}
The proof of Lemma \ref{LEMMA-1-1} is finished.\hfill\qed

Next, we consider the following predictable cumulant process $\Psi (\lambda )=(\Psi
_k(\lambda ),\mathcal{F}_k)_{k=0,...,n}$:
\begin{equation}
\Psi _k ( \lambda )=\sum_{i=1}^k\log \mathbf{E}\big[ e^{\lambda \xi _i}|\mathcal{F}_{i-1} \big]. \label{f29}
\end{equation}
The following lemma gives   a two-sided bound for the process $\Psi(\lambda ).$
\begin{lemma}
\label{LEMMA-1-2} Assume conditions (A1) and (A2).  Then for all $0\leq\lambda\leq   \varepsilon_n^{-1},$
\begin{equation}
\left| \Psi _n(\lambda )-\frac{\lambda ^2}2\right| \leq  c  \,  \lambda^{2+\rho} \varepsilon_n^{ \rho}   +\frac{\lambda^2 \delta_n^2}{2}.
\end{equation}
\end{lemma}

\noindent\textbf{Proof.}
Using a two-term  Taylor's expansion of $\log(1+x), x\geq0$,  we have
\begin{eqnarray*}
 \Psi _n(\lambda ) - \frac{\lambda ^2}2 \langle X \rangle_n
&=& \sum_{i=1}^n\Big( \mathbf{E}[e^{\lambda \xi _i}|\mathcal{F}%
_{i-1}]-1-\lambda \mathbf{E}[\xi _i|\mathcal{F}_{i-1}]-\frac{\lambda ^2}2\mathbf{E}[\xi _i^2|
\mathcal{F}_{i-1}] \Big) \\
&& - \sum_{i=1}^n\frac1{2 \, \big(1+ \theta_i \left(\mathbf{E}[e^{\lambda \xi _i}|\mathcal{F}%
_{i-1}]-1 \right) \big)^2} \left(  \mathbf{E}[e^{\lambda \xi _i}|\mathcal{F}
_{i-1}]-1 \right)^2,
\end{eqnarray*}
where $\theta_i \in (0, 1).$
Since $\mathbf{E}[\xi _i|\mathcal{F}_{i-1}]=0$ and $\mathbf{E}[e^{\lambda \xi _i}|\mathcal{F}_{i-1}]\geq1$ for all $0\leq\lambda\leq     \varepsilon_n^{-1} $,   we deduce that  for  all $0\leq\lambda\leq     \varepsilon_n^{-1}  $,
\begin{eqnarray*}
 \left|\Psi _n(\lambda ) - \frac{\lambda ^2}2 \langle X \rangle _n \right|
&\leq& \sum_{i=1}^n\left| \mathbf{E}[e^{\lambda \xi _i}|\mathcal{F}%
_{i-1}]-1-\lambda \mathbf{E}[\xi _i|\mathcal{F}_{i-1}]-\frac{\lambda ^2}2\mathbf{E}[\xi _i^2|%
\mathcal{F}_{i-1}] \right|\\
&& + \frac12 \sum_{i=1}^n \left(  \mathbf{E}[e^{\lambda \xi _i}|\mathcal{F}%
_{i-1}]-1 \right)^2 .
\end{eqnarray*}
Using condition  (A1) and the inequalities (\ref{inds35})-(\ref{sbopdf}), we get for all $0\leq\lambda\leq   \varepsilon_n^{-1},$
\begin{eqnarray*}
 \left|\Psi _n(\lambda ) - \frac{\lambda ^2}2 \langle X \rangle _n \right|
&\leq&   \sum_{i=1}^{n}\mathbf{E}[ e^{ \lambda \xi_i^+} |\lambda \xi_{i}|^{2+\rho} |\mathcal{F}_{i-1}]
 + \frac12 \sum_{i=1}^n \left(  \mathbf{E}[e^{\lambda \xi _i}|\mathcal{F}%
_{i-1}]-1 \right)^2\\
&\leq&  \lambda^{2+\rho}  \varepsilon_n^{ \rho}   \sum_{i=1}^{n}\mathbf{E}[ \xi_{i} ^ 2 |\mathcal{F}_{i-1}] +    (\lambda\varepsilon_n )^2    \sum_{i=1}^{n}\left(  \mathbf{E}[e^{\lambda \xi _i}|\mathcal{F}%
_{i-1}]-1 \right)\\
&\leq&  \lambda^{2+\rho}  \varepsilon_n^{ \rho}  \langle X\rangle_n + c_1   \lambda^4  \varepsilon_n ^2  \langle X\rangle_n.
\end{eqnarray*}
 Thus
\[
\left|\Psi _n(\lambda ) - \frac{\lambda ^2}2 \langle X \rangle _n \right|  \leq  \Big(1+ c_1(\lambda\varepsilon_n )^{2-\rho} \Big) \lambda^{2+\rho} \varepsilon_n^{ \rho}   \langle X\rangle_n.
\]
Combining the last inequality with condition (A2) and the fact $\langle X \rangle_{n}  \leq 2$,   we get for all $0\leq\lambda\leq   \varepsilon_n^{-1},$
\[
\left|\Psi _n(\lambda ) - \frac{\lambda ^2}{2} \right| \leq
2 \, \Big(1+c_1(\lambda\varepsilon_n )^{2-\rho} \Big) \lambda^{2+\rho} \varepsilon_n^{ \rho}  + \frac{\lambda^2 \delta_n^2}2,
\]
which completes the proof of Lemma \ref{LEMMA-1-2}.\hfill\qed

In the proof of Theorem \ref{co0}, we make use of the following lemma, which gives
us some rates of convergence in the central limit theorem for the conjugate
martingale $Y( \lambda  )$ under the probability measure $\mathbf{P}_{ \lambda  }.$
\begin{lemma}
\label{LEMMA4}
Assume  conditions (A1) and (A2).
\begin{itemize}
   \item[\emph{[i]}] If $\rho \in (0, 1)$, then  there is a positive constant $\alpha$  such that for all $0 \leq \lambda \leq \alpha\,  \varepsilon_n^{-1}  ,$
\[
\sup_{x}\Big| \mathbf{P}_\lambda (  Y_n(\lambda )\leq x)-\Phi (x)\Big| \leq
c_{\alpha, \rho} \, \Big(   ( \lambda \varepsilon_n)^\rho +\varepsilon_n^\rho   +\delta_n \Big) .
\]
In particular, it implies that
\begin{eqnarray} \label{ctlr}
\sup_{x}\Big| \mathbf{P} (  X_n \leq x)-\Phi (x)\Big| \leq
c_{\alpha, \rho} \, \Big(   \varepsilon_n^\rho   +\delta_n \Big) .
\end{eqnarray}
   \item[\emph{[ii]}] If $\rho =1$, then  there is a positive constant $\alpha$  such that for all $0 \leq \lambda \leq \alpha\,  \varepsilon_n^{-1}  ,$
\[
\sup_{x}\Big| \mathbf{P}_\lambda (  Y_n(\lambda )\leq x)-\Phi (x)\Big| \leq
c_{\alpha} \, \Big(    \lambda \varepsilon_n +\varepsilon_n  |\ln \varepsilon_n| +\delta_n \Big) .
\]
In particular, it implies that
\begin{eqnarray} \label{ctlr2}
\sup_{x}\Big| \mathbf{P} (  X_n \leq x)-\Phi (x)\Big| \leq
c_{\alpha} \, \Big(   \varepsilon_n |\log \varepsilon_n |   +\delta_n \Big) .
\end{eqnarray}
 \end{itemize}
\end{lemma}

 The proof of Lemma
\ref{LEMMA4} is complicated, and it is a refinement of the proof of Lemma 3.1 in \cite{FGL13}.
Thus we give details in the supplemental article \cite{FGLSs17}.

\section{Proof of Theorem  \ref{co0}}\label{sec3.0}
\setcounter{equation}{0}
Theorem \ref{co0} will be deduced by the combination of the following two propositions (1 and 2), which are stated
and proved respectively in Subsections  \ref{subsec51} and \ref{sec3.2}. The proof of
the propositions are similar to the proofs  of Theorems 2.1 and 2.2  of  Fan et al.\ \cite{FGL13}.
However, Fan et al.\ \cite{FGL13}  considered the particular case where $\rho=1 .$

\subsection{Proof of  upper bound in Theorem  \ref{co0}}\label{subsec51}
The following   assertion  gives an upper  bound for moderate deviation probabilities.
\begin{proposition}\label{th0}
Assume conditions (A1) and (A2).
\begin{itemize}
   \item[\emph{[i]}] If $\rho \in (0, 1)$, then there is a constant $\alpha >0$  such that
for all $0\leq x \leq \alpha  \varepsilon_n^{-1},$
\begin{equation}\label{t0ie1}
 \frac{\mathbf{P}(X_n>x)}{1-\Phi \left( x\right)}\leq \exp \bigg\{ c_{\alpha,\rho,1} \Big( x^{2+\rho} \varepsilon_n^\rho + x^2 \delta_n^2+(1+ x)\left( \varepsilon_n^\rho   + \delta_n \right)  \Big)  \bigg\} .
\end{equation}
   \item[\emph{[ii]}] If $\rho=1$,   then there is a constant $\alpha >0$  such that
for all $0\leq x \leq \alpha  \varepsilon_n^{-1},$
\begin{equation}\label{t0ie2}
 \frac{\mathbf{P}(X_n>x)}{1-\Phi \left( x\right)}\leq \exp \bigg\{ c_{\alpha,1,1}  \Big( x^{3} \varepsilon_n  + x^2 \delta_n^2 +(1+ x)\left( \varepsilon_n |\ln \varepsilon_n|  + \delta_n \right)\Big)  \bigg\} .
\end{equation}
 \end{itemize}
\end{proposition}

\noindent\textbf{Proof.} For all $0\leq x < 1$, the assertion follows from  (\ref{ctlr}) and (\ref{ctlr2}).
It remains to prove Proposition \ref{th0} for all $1 \leq  x \leq \alpha \, \varepsilon_n^{-1}.$
Changing the probability measure according to (\ref{f21}),  we get for all $ 0\leq\lambda\leq   \varepsilon_n^{-1},$
\begin{eqnarray}
\mathbf{P}(X_n>x) &= & \mathbf{E}_\lambda \left[ Z_n (\lambda)^{-1}\mathbf{1}_{\{ X_n>x \}} \right] \nonumber\\
&= & \mathbf{E}_\lambda \left[\exp \left\{-\lambda X_n+\Psi _n(\lambda )\right\} \mathbf{1}_{\{ X_n>x \}} \right] \nonumber\\
&= & \mathbf{E}_\lambda \left[\exp \left\{
-\lambda Y_n(\lambda)-\lambda B_{n}(\lambda)+\Psi _n(\lambda)\right\} \mathbf{1}_{\{Y_n(\lambda)+B_{n}(\lambda)>x\}} \right].\label{f32}
\end{eqnarray}
Let
$\overline{\lambda}=\overline{\lambda}(x)$ be the
positive solution of the following equation
\begin{equation}\label{f33}
\lambda +\lambda \delta_n^{2} +c  \lambda^{1+\rho} \varepsilon_n^{ \rho}  =x,
\end{equation}
where $c $ is given by inequality (\ref{f25sa}).
The definition of $\overline{\lambda}$ implies that there exist $c_{\alpha,0}, c_{\alpha,1} >0,$ such that for all $1 \leq x \leq \alpha\, \varepsilon_n^{-1},$
\begin{equation}\label{f34}
c_{\alpha,0}\,x \leq \overline{\lambda}  \leq  x
\end{equation}
and
\begin{equation}\label{f35}
\overline{\lambda}=x - c_{\alpha,1}|\theta|(  x^{1+\rho} \varepsilon_n^{ \rho}  +x \delta_n ^2) \in [c_{\alpha,0}, \, \alpha\, \varepsilon_n^{-1}\,].
\end{equation}
By Lemma \ref{LEMMA-1-1}, it follows that $B_{n}(\overline{\lambda})\leq x$.
 From (\ref{f32}), by Lemma \ref{LEMMA-1-2} and equality (\ref{f33}),
 we deduce that for all $1 \leq x \leq \alpha \, \varepsilon_n^{-1},$
\begin{eqnarray}
\mathbf{P}(X_n>x)\leq  e^{ c_{\alpha,2} \,(\overline{\lambda}^{2+\rho}\varepsilon_n^\rho+\overline{\lambda}^2\delta_n^2) -\overline{\lambda}^2/2}\mathbf{E}_{\overline{\lambda}}\big[e^{-
\overline{\lambda}Y_n(\overline{\lambda})}\mathbf{1}_{\{ Y_n(\overline{\lambda})>0\}} \big]. \label{gsfsq}
\end{eqnarray}
Clearly, it holds
\begin{equation}\label{f37d1}
\mathbf{E}_{\overline{\lambda}} \big[ e^{-%
\overline{\lambda}Y_n(\overline{\lambda})}\mathbf{1}_{\{ Y_n(\overline{\lambda})>0\}}\big]= \int_{0}^{\infty} \overline{\lambda} e^{-\overline{\lambda} y}  \mathbf{P}_{\overline{\lambda}}(0 < Y_n(\overline{\lambda})\leq y  ) dy.
\end{equation}
Similarly, for a standard normal random variable $\mathcal{N}$, we have
\begin{equation}\label{f37d2}
\mathbf{E} \big[ e^{-%
\overline{\lambda}\mathcal{N}}\mathbf{1}_{\{ \mathcal{N}>0\}} \big]= \int_{0}^{\infty} \overline{\lambda} e^{-\overline{\lambda} y}   \mathbf{P} (0 < \mathcal{N} \leq y  ) dy.
\end{equation}
From (\ref{f37d1}) and (\ref{f37d2}), it is easy to see that
\begin{eqnarray}
\bigg|\mathbf{E}_{\overline{\lambda}}\big[e^{-%
\overline{\lambda}Y_n(\overline{\lambda})}\mathbf{1}_{\{ Y_n(\overline{\lambda})>0\}}\big]- \mathbf{E} \big[ e^{-%
\overline{\lambda}\mathcal{N}}\mathbf{1}_{\{ \mathcal{N}>0\}}\big] \bigg| \leq
2\sup_y \bigg| \mathbf{P}_{\overline{\lambda}} (Y_n(\overline{\lambda} )\leq y)-\Phi (y) \bigg|.\nonumber
\end{eqnarray}
Using Lemma \ref{LEMMA4}, we get the following bound: for all $
1 \leq  x  \leq \alpha \, \varepsilon_n^{-1},$
\begin{equation}\label{f38}
\bigg|\mathbf{E}_{\overline{\lambda}}\big[e^{-%
\overline{\lambda}Y_n(\overline{\lambda})}\mathbf{1}_{\{ Y_n(\overline{\lambda})>0\}}\big]- \mathbf{E}\big[ e^{-%
\overline{\lambda}\mathcal{N}}\mathbf{1}_{\{ \mathcal{N}>0\}}\big]\bigg| \leq  c_\rho \Big( (\overline{\lambda}\varepsilon_n)^\rho +\widetilde{\varepsilon}_n +\delta_n \Big),
\end{equation}
where
\begin{eqnarray}\label{epsilons}
\widetilde{\varepsilon}_n = \left\{ \begin{array}{ll}
\varepsilon_n^\rho, & \textrm{\ \ \ if $\rho \in (0, 1),$}\\
\varepsilon_n |\ln \varepsilon_n|, & \textrm{\ \ \ if $\rho = 1.$}
\end{array} \right.
\end{eqnarray}
From (\ref{gsfsq}) and (\ref{f38}), we deduce that for all $
1 \leq  x  \leq \alpha \, \varepsilon_n^{-1}, $
\[
\mathbf{P}(X_n>x)\leq e^{ c_{\alpha,2} \,(\overline{\lambda}^{2+\rho}\varepsilon_n^\rho+\overline{\lambda}^2\delta_n^2) -\overline{\lambda}^2/2}\bigg( \mathbf{E}  \big[ e^{- \overline{\lambda}\mathcal{N}}\mathbf{1}_{\{ \mathcal{N}>0\}}  \big]+c_\rho \Big( (\overline{\lambda}\varepsilon_n)^\rho +\widetilde{\varepsilon}_n +\delta_n \Big)  \bigg).
\]
Since
\begin{equation}\label{fsphi}
e^{- \lambda^2/2}\mathbf{E} \left[ e^{-
 \lambda\mathcal{N}}\mathbf{1}_{\{ \mathcal{N}>0\}} \right] =\frac{1}{\sqrt{2\pi}}\int_0^{\infty}e^{-(y+\lambda)^2/2}
  dy=1-\Phi \left(  \lambda\right)
\end{equation}
and
\begin{equation}
1-\Phi \left(  \lambda\right)
\geq
\frac 1{\sqrt{2 \pi}(1+ \lambda)}\ e^{- \lambda^2/2}
\geq
\frac{c_{\alpha,0}}{\sqrt{2 \pi}(1+ c_{\alpha,0})} \frac 1{ \lambda }e^{- \lambda^2/2}, \ \ \ \lambda\geq c_{\alpha,0},
\label{f39}
\end{equation}
we have the following upper bound   for moderate deviation probabilities:
for all $1\leq x \leq \alpha \, \varepsilon_n^{-1} ,$
\begin{eqnarray}
\frac{\mathbf{P}(X_n>x)}{1-\Phi \left( \overline{\lambda}\right)}
&\leq& e^{ c_{\alpha,2} \,(\overline{\lambda}^{2+\rho}\varepsilon_n^\rho+\overline{\lambda}^2\delta_n^2) } \left(\, 1+ c_{\alpha,\rho, 3}\, (    \overline{\lambda}^{1+\rho} \varepsilon_n^\rho +\overline{\lambda}\widetilde{\varepsilon}_n +\overline{\lambda}\delta_n \, ) \right).  \label{f40}
\end{eqnarray}
Next, we would like to make a comparison between $1-\Phi (\overline{\lambda})$ and $1-\Phi (x)$.
  By (\ref{f34}), (\ref{f35}) and (\ref{f39}),  it follows  that
\begin{eqnarray}
   1 \leq \frac{\int_{\overline{\lambda}}^{ \infty}\exp\{- t^2/2 \}d t}{\int_{x}^{ \infty}\exp\{- t^2/2 \}d t}   &= &
  1+\frac{\int_{\overline{\lambda}}^{x}\exp\{ -t^2/2 \} d t}{\int_{x}^{ \infty}\exp\{-t^2/2\}d t}\nonumber\\
   & \leq & 1+c_{\alpha,4}x(x-\overline{\lambda}) \exp\{ (x^2-\overline{\lambda}^2)/2 \}\nonumber\\
   & \leq & \exp\{ c_{\alpha,5}\, (x^{2+\rho}\varepsilon_n^\rho + x^2 \delta_n^2)\}.\label{f41}
\end{eqnarray}
So, it holds
\begin{equation}\label{f42}
1-\Phi \left( \overline{\lambda}\right) =\big( 1-\Phi (x)\big)\exp \left\{ |\theta_{1}| c_{\alpha,5}\, ( x
^{2+\rho}\varepsilon_n^\rho+ x^2 \delta_n^2 ) \right\}.
\end{equation}
Implementing (\ref{f42}) in (\ref{f40}) and using (\ref
{f34}), we obtain for all $1\leq x \leq \alpha \, \varepsilon_n^{-1}, $
\begin{eqnarray*}
\frac{\mathbf{P}(X_n>x)}{1-\Phi \left( x\right) }&\leq& \exp\Big\{ c_{\alpha, 6} ( x^{2+\rho}\varepsilon_n^\rho + x^2 \delta_n^2  ) \Big \}
\left(\frac{}{} 1+c_{\alpha,\rho,7}\left( x^{1+\rho}\varepsilon_n^\rho +x \widetilde{\varepsilon}_n +x\delta_n \right) \right)\\
&\leq & \exp \Big\{ c_{\alpha,\rho,8} \Big( x^{2+\rho}\varepsilon_n^\rho + x^2 \delta_n^2 + x\left( \widetilde{\varepsilon}_n+ \delta_n \right) \Big) \Big \} .
\end{eqnarray*}
This  completes the proof of  Proposition \ref{th0}. \hfill\qed

\subsection{Proof of lower bound in Theorem   \ref{co0}}\label{sec3.2}
The following   assertion  gives a lower  bound for moderate deviation probabilities.
\begin{proposition}\label{th1}
Assume conditions (A1) and (A2).
\begin{itemize}
   \item[\emph{[i]}] If $\rho \in (0, 1)$, then there is a constant $\alpha >0$  such that
for all $0\leq x \leq \alpha  \varepsilon_n^{-1},$
\begin{equation} \label{t1ie1}
\frac{\mathbf{P}(X_n>x)}{1-\Phi \left( x\right)}\geq \exp\bigg\{-c_{\alpha,\rho,2} \,  \Big(  x^{2+\rho} \varepsilon_n^\rho   + x^2 \delta_n^2 + (1+ x) \left( \varepsilon_n^\rho  +  \delta_n \right)\Big) \bigg\}.
\end{equation}
   \item[\emph{[ii]}] If $\rho =1$, then there is a constant $\alpha >0 $ such that
for all $0\leq x \leq \alpha  \varepsilon_n^{-1}$,
\begin{equation}\label{t1ie2}
\frac{\mathbf{P}(X_n>x)}{1-\Phi \left( x\right)}\geq \exp\bigg\{-c_{\alpha,1,2}  \Big(  x^{3} \varepsilon_n   + x^2 \delta_n^2 + (1+ x) \left( \varepsilon_n |\ln \varepsilon_n| +  \delta_n \right)\Big) \bigg\}.
\end{equation}
 \end{itemize}
\end{proposition}

\noindent\textbf{Proof.}
For all $0\leq x < 1$, the assertion follows from (\ref{ctlr}) and (\ref{ctlr2}).
It remains to prove Proposition \ref{th1} for all $1 \leq x \leq \alpha  \varepsilon_n^{-1}$, where $\alpha >0$ is  a small constant.
 Let
$\underline{\lambda}=\underline{\lambda}(x)$ be the smallest positive solution of the following equation
\begin{equation}\label{f44}
\lambda -\lambda \delta_n^{2} -c  \lambda^{1+\rho}\varepsilon_n^\rho  =x,
\end{equation}
where   $c$ is given   by inequality (\ref{f25sa}). 
The definition of $\underline{\lambda}$ implies that for all $1 \leq x \leq  \alpha \varepsilon_n^{-1}, $
\begin{equation}\label{f45}
x\leq \underline{\lambda}\leq   c_{\alpha,1}\,x
\end{equation}
and
\begin{equation}\label{f46}
\underline{\lambda}=x+c_{\alpha,2}|\theta|(x^{1+\rho}\varepsilon_n^\rho +x \delta_n ^2) \in [1, \      \varepsilon_n^{-1}].
\end{equation}
 From (\ref{f32}), using Lemmas \ref{LEMMA-1-1},  \ref{LEMMA-1-2} and equality (\ref{f44}),
  we have for all $1 \leq x \leq   \alpha \varepsilon_n^{-1}, $
\begin{equation}\label{jknjssa}
\mathbf{P}(X_n>x)\geq e^{ - c_{1} \,(\underline{\lambda}^{2+\rho}\varepsilon_n^\rho+\underline{\lambda}^2\delta_n^2) -\underline{\lambda}^2/2} \mathbf{E}_{\underline{\lambda}} \left[e^{-%
\underline{\lambda}Y_n(\underline{\lambda})}\mathbf{1}_{\{ Y_n(\underline{\lambda})>0 \}} \right].
\end{equation}

In the subsequent we distinguish $\underline{\lambda}$ into two cases.
First, let $1\leq \underline{\lambda}\leq \alpha_1 \min\{ \varepsilon_n^{-\rho/(1+\rho)}  , \delta_n^{-1}  \}$, where $\alpha_1>0$ is a small positive constant
whose exact value will be given later.
Note that inequality (\ref{f38}) can be established with $\overline{\lambda}$ replaced by $\underline{\lambda}$, which, in turn,
implies that
\[
\mathbf{P}(X_n>x)\geq e^{ - c_{1} \,(\underline{\lambda}^{2+\rho}\varepsilon_n^\rho+\underline{\lambda}^2\delta_n^2) -\underline{\lambda}^2/2} \bigg( \mathbf{E} \left[e^{-\underline{\lambda}  \mathcal{N}}\mathbf{1}_{\{ \mathcal{N}>0\}} \right]-c_{\rho,2}\Big(  (\underline{\lambda}\varepsilon_n)^\rho +\widetilde{\varepsilon}_n +\delta_n \Big)  \bigg),
\]
where $\widetilde{\varepsilon}_n$ is defined  by  (\ref{epsilons}).
By (\ref{fsphi}) and (\ref{f39}),
we get the following lower bound on tail probabilities:
\begin{eqnarray}
\frac{\mathbf{P}(X_n>x)}{1-\Phi \left( \underline{\lambda}\right)}
&\geq& e^{ - c_{1} \,(\underline{\lambda}^{2+\rho}\varepsilon_n^\rho+\underline{\lambda}^2\delta_n^2) } \left(\frac{}{} 1-c_{\rho,2}\left(  \underline{\lambda}^{1+\rho} \varepsilon_n^\rho +\underline{\lambda} \widetilde{\varepsilon}_n +\underline{\lambda} \delta_n^2    \right) \right).  \label{f51}
\end{eqnarray}
Taking  $\alpha_1 = \min\{ \frac{1}{ (8c_{\rho,2})^{ 1/(1+\rho)}},   \frac{1}{8c_{\rho,2}}     \}$,
we have for all $1\leq \underline{\lambda} \leq \alpha_1 \min \{\varepsilon_n^{-\rho/(1+\rho)}, \delta_n^{-1} \} $,
\begin{eqnarray}\label{f55f}
1-c_{\rho,2}\left( \underline{\lambda}^{1+\rho} \varepsilon_n^\rho +\underline{\lambda}\widetilde{\varepsilon}_n +\underline{\lambda}\delta_n \right)&\geq &  \exp\left\{-2c_{\rho,2}\left( \underline{\lambda}^{1+\rho} \varepsilon_n^\rho +\underline{\lambda}\widetilde{\varepsilon}_n +\underline{\lambda}\delta_n \right) \right\} .
\end{eqnarray}
Implementing (\ref{f55f}) in (\ref{f51}), we get
\begin{eqnarray}
\frac{\mathbf{P}(X_n>x)}{1-\Phi \left( \underline{\lambda}\right)}&\geq& \exp \bigg \{ -c_{\rho,3} \left( \underline{\lambda}^{2+\rho}\varepsilon_n^\rho +\underline{\lambda }\widetilde{\varepsilon}_n +\underline{\lambda}\delta_n +\underline{\lambda}^2\delta_n^2 \right)\bigg\}   \label{f54}
\end{eqnarray}
which holds for all $1\leq \underline{\lambda} \leq \alpha_1 \min \{ \varepsilon_n^{-\rho/(1+\rho)}, \delta_n^{-1}\} $.

Next,  consider the case of $\alpha_1 \min\{ \varepsilon_n^{-\rho/(1+\rho)} , \delta_n^{-1} \} \leq \underline{\lambda} \leq \alpha \varepsilon_n ^{-1}.$ Let $K \geq 1$ be a constant, whose exact value will be chosen later.
It is obvious that
\begin{eqnarray}\label{jknjsta}
\mathbf{E}_{\underline{\lambda}} \left[e^{-\underline{\lambda}Y_n(\underline{\lambda})}\mathbf{1}_{\{ Y_n(\underline{\lambda})>0 \}} \right] &\geq& \mathbf{E}_{\underline{\lambda}} \Big[e^{-\underline{\lambda}Y_n(\underline{\lambda})}\mathbf{1}_{\{0< Y_n(\underline{\lambda})\leq  K \tau \}} \Big] \nonumber\\
 &\geq&e^{-\underline{\lambda} K \tau}\mathbf{P}_{\underline{\lambda}} \Big(0< Y_n(\underline{\lambda})\leq  K \tau \Big),
\end{eqnarray}
where $\tau = (\underline{\lambda}\varepsilon_n)^\rho +  \widetilde{\varepsilon}_n  +\delta_n$.
 From Lemma \ref{LEMMA4}, we get
\begin{eqnarray*}
\mathbf{P}_{\underline{\lambda}} \Big(0< Y_n(\underline{\lambda})\leq  K \tau \Big) &\geq&   \mathbf{P}  \Big( 0<  \mathcal{N}
\leq  K \tau  \Big)  -  c_{\rho,5} \tau  \\
  &\geq& K \tau  e^{-K^2  \tau^2/2}   - c_{\rho,5} \tau\\
  &\geq& \left( K e^{-8 K^2  \alpha}   - c_{\rho,5}\right) \tau.
\end{eqnarray*}
Taking $\alpha= 1/(16 K^2) $, we obtain
 $$\mathbf{P}_{\underline{\lambda }} \Big(0< Y_n(\underline{\lambda})\leq  K \tau \Big)  \geq     \left( \frac12 K   - c_{\rho,5} \right) \tau.$$
Letting $K\geq   8c_{\rho, 5}$, we deduce that
$$\mathbf{P}_{\underline{\lambda}} \Big(0< Y_n(\underline{\lambda })\leq  K \tau \Big)  \geq \frac38 K \tau \geq \frac38 K  \frac{  \max \left\{ \underline{\lambda }^{1+\rho}\varepsilon_n^\rho  , \underline{\lambda } \delta_n  \right\}} { \underline{\lambda } }.$$
Choosing $K=\max \Big\{8 c_{\rho,5},  \frac{16\alpha_1^{-1-\rho}}{3\sqrt{\pi}} \Big\}$ and taking into account that
$ \alpha_1 \min\{ \varepsilon_n ^{-\rho/(1+\rho)},  \delta_n^{-1} \} \leq \underline{\lambda} \leq \alpha \varepsilon_n^{-1}$,
we get
\begin{eqnarray*}
\mathbf{P}_{\underline{\lambda }} \Big(0< Y_n(\underline{\lambda }) \leq  K \tau \Big)  \geq    \frac{2}{\sqrt{\pi}\underline{\lambda} } .
\label{jknjstb}
\end{eqnarray*}
Since the inequality
$$\frac{2}{\sqrt{\pi}\lambda }  e^{- \lambda^2/2}  \geq 1-\Phi \left(  \lambda \right)$$ is valid for all $\lambda >0$,
it follows that for all $ \alpha_1 \min\{ \varepsilon_n ^{-\rho/(1+\rho)},  \delta_n^{-1} \} \leq \underline{\lambda} \leq \alpha \varepsilon_n^{-1}$,
\begin{eqnarray}\label{dfac}
\mathbf{P}_{\underline{\lambda}} \Big(0< Y_n(\underline{\lambda})\leq  K \tau\Big)   \geq \Big( 1-\Phi \left(  \underline{\lambda} \right) \Big)e^{\underline{\lambda}^2/2} .
\end{eqnarray}
From (\ref{jknjssa}), (\ref{jknjsta}) and (\ref{dfac}), we get
 \begin{eqnarray}
\frac{\mathbf{P}(X_n>x)}{1-\Phi \left( \underline{\lambda}\right)}&\geq& \exp \bigg \{ -c_{\alpha, 6} \left( \underline{\lambda}^{2+\rho}\varepsilon_n^\rho +\underline{\lambda }\widetilde{\varepsilon}_n +\underline{\lambda}\delta_n +\underline{\lambda}^2\delta_n^2\right)\bigg\}  \label{fgj53}
\end{eqnarray}
which holds for all $ \alpha_1 \min\{ \varepsilon_n ^{-\rho/(1+\rho)},  \delta_n^{-1} \} \leq \underline{\lambda} \leq \alpha \varepsilon_n^{-1}$.

Combining (\ref{f54}) and (\ref{fgj53}) together, we obtain for all $1\leq \underline{\lambda} \leq  \alpha \varepsilon_n^{-1}     ,$
\begin{eqnarray}\label{ft52}
\frac{\mathbf{P}(X_n>x)}{1-\Phi \left( \underline{\lambda}\right)}&\geq& \exp \bigg \{ -c_{\alpha,\rho, 7} \left( \underline{\lambda}^{2+\rho}\varepsilon_n^\rho +\underline{\lambda }\widetilde{\varepsilon}_n +\underline{\lambda}\delta_n +\underline{\lambda}^2\delta_n^2\right)\bigg\}  .
\end{eqnarray}
By a similar argument as in (\ref{f41}), it is easy to see that
\begin{equation}\label{f53}
1-\Phi \left( \underline{\lambda}\right) =\Big( 1-\Phi (x)\Big)  \exp \left\{ - |\theta| c_3\, ( x
^{2+\rho}\varepsilon_n^\rho+ x^2 \delta_n^2  ) \right\} .
\end{equation}
Combining (\ref{f45}), (\ref{ft52}) and (\ref{f53}) together, we find that  for all $1\leq x \leq \alpha  \varepsilon_n^{-1}   ,$
 \begin{equation}
\frac{\mathbf{P}(X_n>x)}{1-\Phi \left( x\right)}  \geq  \exp \bigg \{ -c_{\alpha, \rho, 8} \left(x^{2+\rho}\varepsilon_n^\rho +x\widetilde{\varepsilon}_n +x\delta_n +x^2\delta_n^2  \right)\bigg\},
\end{equation}
 which gives the conclusion of Proposition \ref{th1}. \hfill\qed

\section{Proof of Corollary \ref{mdp}} \label{sec3.3}
\setcounter{equation}{0}
To prove Corollary \ref{mdp}, we need the following two-sides bound  on tail probabilities of the standard normal random variable:
\begin{eqnarray}\label{fgsgj1}
\frac{1}{\sqrt{2 \pi}(1+x)} e^{-x^2/2} \leq 1-\Phi ( x ) \leq \frac{1}{\sqrt{ \pi}(1+x)} e^{-x^2/2}, \ \ \ \   x\geq 0.
\end{eqnarray}
First, we prove that for any given Borel set $B\subset [0, \infty),$
\begin{eqnarray}\label{dfgkmsf}
 \limsup_{n\rightarrow \infty}\frac{1}{a_n^2}\log \mathbf{P}\left(\frac{1}{a_n} X_n \in B \right) \leq  - \inf_{x \in \overline{B}}\frac{x^2}{2}.
\end{eqnarray}
Let $x_0=\inf_{x\in B} x.$ Then it is obvious that $x_0 \geq 0$ and  $x_0\geq\inf_{x\in \overline{B}} x.$ By Theorem \ref{co0}, we deduce that
\begin{eqnarray*}
&&\mathbf{P}\left(\frac{1}{a_n} X_n \in B \right)\\
 &&\leq  \mathbf{P}\left( X_n  \geq a_n x_0\right)\\
 &&\leq  \Big( 1-\Phi \left( a_nx_0\right)\Big)\exp\bigg\{c_{\alpha} \Big(  \left( a_nx_0\right)^{2+\rho} \varepsilon_n^\rho   + \left( a_nx_0\right)^2 \delta_n^2 + (1+ \left( a_nx_0\right)) \left( \varepsilon_n^\rho|\ln \varepsilon_n | +  \delta_n \right)\Big) \bigg\}.
\end{eqnarray*}
Using  (\ref{fgsgj1}) and the assumption $a_n\epsilon_n \rightarrow 0,$
we have
\begin{eqnarray*}
\limsup_{n\rightarrow \infty}\frac{1}{a_n^2}\log\mathbf{P}\left(\frac{1}{a_n} X_n \in B \right)
 \ \leq \  -\frac{x_0^2}{2} \ \leq \  - \inf_{x \in \overline{B}}\frac{x^2}{2} ,
\end{eqnarray*}
which gives (\ref{dfgkmsf}).

Next, we prove that for any given Borel set $B\subset [0, \infty),$
\begin{eqnarray}\label{dfgk02}
\liminf_{n\rightarrow \infty}\frac{1}{a_n^2}\log \mathbf{P}\left(\frac{1}{a_n} X_n \in B \right) \geq   - \inf_{x \in B^o}\frac{x^2}{2} .
\end{eqnarray}
For any $\varepsilon_1>0,$ there exists an $x_0 \in B^o,$ such that
\begin{eqnarray}
 \frac{x_0^2}{2} \leq   \inf_{x \in B^o}\frac{x^2}{2} +\varepsilon_1.
\end{eqnarray}
For $x_0 \in B^o,$ there exists an $\varepsilon_2 > 0,$ such that $(x_0-\varepsilon_2, x_0+\varepsilon_2]  \subset B.$
Then it is obvious that $x_0\geq\inf_{x\in  B^o } x.$
By Theorem \ref{co0}, we deduce that
\begin{eqnarray*}
&&\mathbf{P}\left(\frac{1}{a_n} X_n \in B \right)
\ \geq\  \mathbf{P}\Big( X_n  \in (a_n ( x_0-\varepsilon_2), a_n( x_0+\varepsilon_2)] \Big)\\
&&   \geq \mathbf{P}\Big( X_n  > a_n ( x_0-\varepsilon_2)   \Big)-\mathbf{P}\Big( X_n >   a_n( x_0+\varepsilon_2) \Big)\\
&&  \geq \Big( 1-\Phi \left( a_n( x_0-\varepsilon_2)\right)\Big)\exp\bigg\{-c_{\alpha} \Big(  \left( a_n( x_0-\varepsilon_2)\right)^{2+\rho} \varepsilon_n^\rho + \left(a_n( x_0-\varepsilon_2)\right)^2 \delta_n^2 \\
&&  \ \ \  \ \ \   \ \ \    \ \ \    \ \ \   \ \ \     \ \ \    \ \ \     \ \ \     \ \ \     \ \ \      \ \ \     \      \ \ \    \ \ \    \ \ \   \ \ \     \ \ \     + (1+ \left( a_n( x_0-\varepsilon_2)\right)) \left( \varepsilon_n^\rho|\ln \varepsilon_n |  +  \delta_n \right)\Big) \bigg\}\\
&&  \ \ \ \   -\Big( 1-\Phi \left( a_n( x_0+\varepsilon_2)\right)\Big)\exp\bigg\{c_{\alpha} \Big(  \left( a_n( x_0+\varepsilon_2)\right)^{2+\rho} \varepsilon_n^\rho  + \left(a_n( x_0+\varepsilon_2)\right)^2 \delta_n^2\\
&& \ \ \   \ \ \  \ \ \   \ \ \  \ \ \ \    \ \ \ \ \ \ \    \ \ \   \ \ \    \ \ \   \ \ \    \     \ \ \     \      \ \ \    \ \ \    \ \ \   \ \ \   + (1+ \left( a_n( x_0+\varepsilon_2)\right)) \left( \varepsilon_n^\rho|\ln \varepsilon_n |  +  \delta_n \right)\Big) \bigg\} \\
&&=: P_{1, n} - P_{2, n}.
\end{eqnarray*}
Since $a_n\epsilon_n \rightarrow 0$, it is easy to see that $\lim_{n\rightarrow \infty} P_{2, n} /P_{1, n} = 0.$ Thus for $n$ large enough, it holds
\begin{eqnarray*}
&&\mathbf{P}\left(\frac{1}{a_n} X_n \in B \right) \geq \frac12 P_{1, n} .
\end{eqnarray*}
Using (\ref{fgsgj1}) and the assumption $a_n\epsilon_n \rightarrow 0$ again, it follows that
\begin{eqnarray*}
\liminf_{n\rightarrow \infty}\frac{1}{a_n^2}\log \mathbf{P}\left(\frac{1}{a_n} X_n \in B \right)  \geq  -  \frac{1}{2}( x_0-\varepsilon_2)^2 . \label{ffhms}
\end{eqnarray*}
Letting $\varepsilon_2\rightarrow 0,$  we get
\begin{eqnarray*}
\liminf_{n\rightarrow \infty}\frac{1}{a_n^2}\log \mathbf{P}\left(\frac{1}{a_n} X_n \in B \right)  \geq  -  \frac{x_0^2}{2}     \geq     -\inf_{x \in B^o}\frac{x^2}{2} -\varepsilon_1.
\end{eqnarray*}
Since $\varepsilon_1$ can be arbitrary small, we obtain (\ref{dfgk02}).   \hfill\qed

\section{ Proof of Theorem \ref{th3.0}   }
\setcounter{equation}{0}
To prove Theorem \ref{th3.0}, we need the following lemma.
\begin{lemma}\label{lemma9}
Assume the conditions of Theorem \ref{th3.0}. Then  for all $x\geq0,$
\begin{eqnarray}
\mathbf{P}\Big(   |W_n | > x  \Big)   \leq  2 \exp \Bigg\{-  \frac{x^2}{2(1+ \frac M  n + \frac{x L}{3\sqrt{n}}  )}  \Bigg\} .
\end{eqnarray}
\end{lemma}
\noindent\textbf{Proof.}
Let $T_0=\min\{ K, \frac L 3 \}.$
It is easy to see that for all $0\leq \lambda < T_0 $,
\begin{eqnarray}
\mathbf{E}[e^{ \lambda \eta_i}| \mathcal{F}_{i-1}] &\leq& 1 +  \lambda \mathbf{E}[\eta_i  | \mathcal{F}_{i-1}] + \frac{\lambda^2}{2} \mathbf{E}[ \eta_i^2 | \mathcal{F}_{i-1}] + \frac{\lambda^3}{3!} \mathbf{E}[ |\eta_i|^3e^{ K |\eta_i|}   | \mathcal{F}_{i-1}]  \nonumber \\
&\leq& 1 +    \frac{\lambda^2}{2}(1+ \frac13\lambda L) \mathbf{E}[ \eta_i^2 | \mathcal{F}_{i-1}]   \nonumber \\
&\leq& \exp\bigg\{   \frac{\lambda^2}{2}(1+ \frac13\lambda L) \mathbf{E}[ \eta_i^2 | \mathcal{F}_{i-1}] \bigg\}    \nonumber \\
&\leq& \exp\bigg\{   \frac{\lambda^2}{2(1- \lambda T_0)}  \mathbf{E}[ \eta_i^2 | \mathcal{F}_{i-1}] \bigg\} , \nonumber
\end{eqnarray}
which implies that for all $0\leq \lambda < T_0, $
\begin{eqnarray*}
&&\mathbf{E}\bigg[\exp\bigg\{ \lambda \sum_{i=1}^n\eta_i -\frac{\lambda^2\, \Xi_n}{2(1- \lambda T_0)} \bigg\}  \bigg]  \\
&&\leq \mathbf{E}\Bigg[\exp\bigg\{ \lambda \sum_{i=1}^{n-1}\eta_i -\frac{\lambda^2\, \Xi_{n-1}}{2(1- \lambda T_0)}   \bigg\}   \mathbf{E}\bigg[ \exp\bigg\{ \lambda  \eta_n -\frac{\lambda^2\mathbf{E}[ \eta_n^2 | \mathcal{F}_{n-1}]}{2(1- \lambda T_0)}  \bigg\} \bigg| \mathcal{F}_{n-1} \bigg]  \Bigg]  \\
&&\leq \mathbf{E}\Bigg[\exp\bigg\{ \lambda \sum_{i=1}^{n-1}\eta_i -\frac{\lambda^2 \, \Xi_{n-1}}{2(1- \lambda T_0)}   \bigg\} \Bigg] \\
&&\leq  1,
\end{eqnarray*}
where
$$ \Xi_n   =\sum_{i=1}^{n}\mathbf{E}[ \eta_i^2 | \mathcal{F}_{i-1}].$$
Since $ \Xi_{n}\leq n+M$ a.s.,  we have for all $x \geq0$ and all $0\leq \lambda < T_0, $
\begin{eqnarray*}
 \mathbf{P}\Big( W_n>x \Big) & =&  \mathbf{P}\Big( \sum_{i=1}^n\eta_i >x \sqrt{n} \Big)  \\
&\leq& \mathbf{E}\bigg[\exp\bigg\{ - \lambda x \sqrt{n} +   \lambda \sum_{i=1}^{n}\eta_i -\frac{\lambda^2\, \Xi_{n}}{2(1- \lambda T_0)}  +\frac{\lambda^2 (n+M)}{2(1- \lambda T_0)}     \bigg\}    \bigg]  \\
&\leq&    \mathbf{E}\bigg[\exp\bigg\{ - \lambda x \sqrt{n}   +\frac{\lambda^2 (n+M)}{2(1- \lambda T_0)}     \bigg\}    \bigg] .
\end{eqnarray*}
Thus for all $x \geq0,$
\begin{eqnarray}
 \mathbf{P}\Big( W_n>x \Big) &\leq&  \inf_{0\leq \lambda < T_0}  \mathbf{E}\bigg[\exp\bigg\{ - \lambda x \sqrt{n}   +\frac{\lambda^2 (n+M)}{2(1- \lambda T_0)}     \bigg\}    \bigg] \nonumber  \\
 &\leq&   \exp\Bigg\{ - \frac{x^2}{2(1+M/n +x T_0/\sqrt{n} )}     \Bigg\} \nonumber \\
  &\leq&   \exp\Bigg\{ - \frac{x^2}{2(1+ \frac M  n + \frac{x L}{3\sqrt{n}} )}     \Bigg\}  .\label{ssfed01}
\end{eqnarray}
Similarly, we have for all $x \geq0,$
\begin{eqnarray}
 \mathbf{P}\Big( W_n < -x \Big)
   \leq    \exp\Bigg\{ - \frac{x^2}{2(1+ \frac M  n + \frac{x L}{3\sqrt{n}} )}     \Bigg\}  . \label{ssfed02}
\end{eqnarray}
Combining (\ref{ssfed01}) and (\ref{ssfed02}) together, we obtain  the desired inequality. \hfill\qed

Now we  are in position to  prove Theorem \ref{th3.0}.
 By Theorem \ref{co0}, there exist constants $\alpha \in (0, 1]$  and $C\geq1$ such that
for all  $ 0 \leq x \leq\alpha \,  n^{1/2} ,$
\begin{eqnarray}\label{gonineq01}
 \frac { \mathbf{P}\left( W_n>x \right) } { 1-\Phi \left( x\right)} = \exp \bigg\{ \theta C (1+x^3) \frac{\ln n}{\sqrt{n}} \bigg\}
\end{eqnarray}
and
\begin{eqnarray}\label{gonineq02}
 \frac { \mathbf{P}\left( W_n<-x \right) } { \Phi \left(- x\right)} = \exp \bigg\{ \theta C (1+x^3) \frac{\ln n}{\sqrt{n}} \bigg\},
\end{eqnarray}
where $|\theta |\leq 1$.
By Theorem 1 of Mason and Zhou \cite{MZ12} with $\varepsilon_n=\alpha$ and $K_n = C \ln n$,  then whenever $n \geq   64 C^2 (\ln n )^2 $ and
$$|W_n | \leq \frac{\sqrt{ n}}{8\ln n},$$
 we have
$$|W_n -Z|\leq   2 C \Big( W_n^2\, +1  \Big) \frac{ \ln n }{\sqrt{n}} ,$$
which gives (\ref{fg235}).  Notice that there exists an integer $n_0$ such that  $n \geq   64 C^2 (\ln n )^2$ for all $n\geq n_0.$

Next we give the proof of (\ref{fg2sfs5}). By (\ref{fg235}), we have for all $0 \leq x \leq \frac{C\,}{32} \, n/(\ln n)^2,  $
\begin{eqnarray}
 \mathbf{P}\Big( \sqrt{n}|W_n-Z|/\ln n > x \Big)&\leq& \mathbf{P}\Big( \sqrt{n}|W_n-Z|/\ln n > x , |W_n | \leq \frac18\, \sqrt{ n}/\ln n\Big) \nonumber\\
 &&+ \, \mathbf{P}\Big(   |W_n | > \frac18\, \sqrt{ n}/\ln n \Big) \nonumber \\
 &\leq& \mathbf{P}\Big( 2 C \big( W_n^2\, +1  \big) > x  \Big) + \, \mathbf{P}\Big(   |W_n | > \frac18\, \sqrt{ n}/\ln n \Big) \nonumber\\
 &\leq&  \mathbf{P}\Big(   |W_n | > \sqrt{ x/(2C) }  \Big) + \, \mathbf{P}\Big(   |W_n | > \frac18\, \sqrt{ n}/\ln n \Big). \label{fsf01}
\end{eqnarray}
Notice that
$$1-\Phi \left( x\right) \leq   \exp\{-x^2/2\}, \ \ x \geq 0.$$
When $ 0 \leq x \leq 2C \alpha^2 n  /(8 C \ln n  )^2, n\geq 2,$ by the inequalities (\ref{gonineq01}) and (\ref{gonineq02}), it holds that
\begin{eqnarray}\label{fsf02}
\mathbf{P}\Big(   |W_n | > \sqrt{ x/(2C) }  \Big)   &\leq&  2 \exp \bigg\{ -\frac14(\sqrt{ x/(2C) })^2  \bigg\} \nonumber\\
&=&    \exp \bigg\{1 -\frac1{8C } x \bigg\},
\end{eqnarray}
and that
\begin{eqnarray}\label{fsf03}
\mathbf{P}\Big(   |W_n | > \frac18\, \sqrt{ n}/\ln n \Big)  &\leq&   2\exp \bigg\{ -\frac n{8\cdot32 (\ln n)^2 }  \bigg\}\nonumber \\
&\leq&   \exp \bigg\{ 1-\frac { C}{8 \alpha^2 } x \bigg\}.
\end{eqnarray}
Returning to (\ref{fsf01}), we obtain  for all $ 0 \leq x \leq 2C \alpha^2 n  /(8 C \ln n   )^2,$
\begin{eqnarray} \label{fsf1032}
\mathbf{P}\Big( \sqrt{n}|W_n-Z|/\ln n > x \Big) &\leq&  2 \exp \Big\{ 1- c' x \Big\} ,
\end{eqnarray}
where $c'= \min\{ \frac1{8C }  , \frac { C}{8 \alpha^2 } \}.$
When $x >2C \alpha^2 n  /(8 C \ln n   )^2,$ it holds
\begin{eqnarray}\label{fssfh023}
 \mathbf{P}\Big( \sqrt{n}|W_n-Z|/\ln n > x \Big) \leq   \mathbf{P}\Big( \sqrt{n}|W_n|/\ln n > x/2 \Big) +  \mathbf{P}\Big( \sqrt{n}|Z|/\ln n > x /2 \Big).
\end{eqnarray}
By Lemma \ref{lemma9},  there exists a positive constant $ \lambda$ such that for all $x >2C \alpha^2 n  /(8 C \ln n   )^2,$
\begin{eqnarray}
\mathbf{P}\Big(   \sqrt{n}|W_n|/\ln n > x/2  \Big)   &\leq&  2 \exp \bigg\{ - \frac{3}{8 L} x  \sqrt{ n}\frac{\ln n }{\sqrt{n}}   \bigg\}  \nonumber  \\
&\leq &   \exp  \bigg\{1- \frac{3}{8 L} x  \bigg\} , \nonumber
\end{eqnarray}
and that
\begin{eqnarray}
\mathbf{P}\Big(   \sqrt{n}|Z|/\ln n > x/2  \Big)   &\leq&  2 \exp \bigg\{- \frac18 x^2 \frac{\ln n }{\sqrt{n}}   \bigg\} \nonumber  \\
&\leq &   \exp  \bigg\{1- \frac{  \alpha^2  }{256 C  } x   \bigg\}.\nonumber
\end{eqnarray}
Returning to (\ref{fssfh023}), we have for all $x >2C \alpha^2 n  /(8 C \ln n   )^2,$
\begin{eqnarray} \label{fghm53}
 \mathbf{P}\Big( \sqrt{n}|W_n-Z|/\ln n > x \Big) \leq  2 \exp  \Big\{ 1- c'' x  \Big\},
\end{eqnarray}
where $c'' = \min\{ \frac{3}{8 L} , \frac{  \alpha^2  }{256 C  } \}.$ Combining (\ref{fsf1032}) and (\ref{fghm53}) together, we get (\ref{fg2sfs5}).

\section{ Proof of Theorem \ref{tghbt}   }
\setcounter{equation}{0}
The main idea of the proof of Theorem \ref{tghbt} is to use $m$-dependence approximation. We make use of the following  lemma of Berbee \cite{B87}.
\begin{lemma}\label{lema7.1}
Let $(Y_i)_{ 1\leq i \leq n} $ be a sequence of random variables on some probability space and define $\beta^{(i)}= \beta(Y_i, (Y_{i+1},...,Y_n) )$.
Then the probability space can be extended with random variables $\widetilde{Y}_i$ distributed as $Y_i$ such  that $(\widetilde{Y}_i)_{ 1\leq i \leq n} $ are independent and
$$\mathbf{P}( Y_i\neq \widetilde{Y}_i\textrm{ for some }1\leq i \leq n ) \leq \beta^{(1)}+ ...+ \beta^{(n-1)}.$$
\end{lemma}

%

Now we are in position to prove Theorem \ref{tghbt}. Recall
   $m=\lfloor n^\alpha \rfloor$ and $k=\lfloor n/(2m) \rfloor.$
By Lemma \ref{lema7.1}, there exists a sequence of independent random variables $(\widetilde{Y}_j)_{1\leq j \leq k}$ such that
$\widetilde{Y}_j$ and $Y_j$ have the same distribution for each $1\leq j \leq k$ and
\begin{equation}\label{ineq3.6}
\mathbf{P}( Y_i\neq \widetilde{Y}_i\textrm{ for some }1\leq i \leq k ) \leq k \beta(m) \leq a_1\exp\{-0.5 a_2 n^{ \alpha \, \tau }\}.
\end{equation}
Therefore, we have
\begin{equation}\label{ineq3.7}
 \big|\mathbf{P}\big(S_n/\sqrt{\mathbf{E}S_n^2} >x \big)- \mathbf{P} \big(\widetilde{S}_n/\sqrt{\mathbf{E}S_n^2} >x \big) \big |  \leq a_1\exp\big\{-0.5 a_2 n^{  \alpha \, \tau } \big\},
\end{equation}
where $ \widetilde{S}_n=\sum_{j=1}^k \widetilde{Y}_j.$
 By (\ref{cosf3.11}) and (\ref{cosdf3.12}), we have
 $$\mathbf{E}  |\widetilde{Y}_{i}|^{2+\rho}   \leq    c_1^{2+\rho}c_2^{-2}  m^{ \rho /2} \mathbf{E}  \widetilde{Y}_{i} ^2   $$
 for all $1\leq j  \leq k, $  and $$ \textrm{Var}  (\widetilde{S}_n)   \asymp n .$$
By   (\ref{ineq3.6}) and (\ref{cosf3.11}), it is easy to see that
\begin{eqnarray} \label{fsss01}
\Big|\mathbf{E} \widetilde{S}_n^2   -  \mathbf{E}S _n^2   \Big |
&=&\Big|\mathbf{E}[ (\widetilde{S}_n^2 -    S _n^2)\mathbf{1}_{\{Y_i\neq \widetilde{Y}_i\textrm{ for some }1\leq i \leq k\}} ] \Big | \nonumber \\
&\leq& 2  \mathbf{E}[ e^{ \frac14 a_2n^{\alpha   \tau}  }  \mathbf{1}_{\{Y_i\neq \widetilde{Y}_i\textrm{ for some }1\leq i \leq k\}}] \nonumber\\
 && + \, \mathbf{E} [\widetilde{S}_n^2 \mathbf{1}_{\{ \widetilde{S}_n > e^{ \frac18a_2n^{\alpha   \tau} }  \}}]    + \mathbf{E} [S_n^2 \mathbf{1}_{\{ S_n > e^{ \frac18a_2n^{\alpha   \tau} }  \}}] \nonumber \\
 &\leq& 2  e^{ \frac14 a_2n^{\alpha   \tau}  } \mathbf{P}(Y_i\neq \widetilde{Y}_i\textrm{ for some }1\leq i \leq k  ) \nonumber\\
 && + \, e^{ -\frac{\rho}8a_2n^{\alpha   \tau} } \mathbf{E}  |\widetilde{S}_n|^{2+\rho}      + e^{ -\frac{\rho}8a_2n^{\alpha   \tau} } \mathbf{E}|S_n|^{2+\rho}  \nonumber \\
    &\leq& O(1)\exp\Big\{-\frac{\rho }{16} a_2 n^{  \alpha   \tau} \Big\}  \nonumber \\
    &=& O( n^{-2}). \nonumber
\end{eqnarray}
It is obvious that $\widetilde{Y}_j \leq n^\alpha c_3$ a.s. Applying  Theorem \ref{co0}   to $\widetilde{S}_n/\sqrt{ \mathbf{E}S _n^2},$ we deduce that
  there is a constant $\alpha >0$, such that
for all $0\leq x =o( n^{\frac12 - \alpha    }),$
\begin{equation}\label{ineq3.8}
\frac{\mathbf{P}( \widetilde{S}_n/\sqrt{\mathbf{E}S _n^2}>x )}{1-\Phi \left( x\right)}=\exp\Bigg\{\theta_1 c_{\rho}  \frac{ (1+x)^{2+\rho }}{n^{\rho (\frac12 - \alpha )} }     \Bigg\}.
\end{equation}
The inequalities (\ref{ineq3.7}) and (\ref{ineq3.8})  together implies that
\begin{eqnarray}
\frac{\mathbf{P}( S_n/\sqrt{\mathbf{E}S _n^2}>x )}{1-\Phi \left( x\right)}&=&\exp\Bigg\{\theta_1 c_{1,\rho}  \frac{ (1+x)^{2+\rho }}{n^{\rho(\frac12 - \alpha)}  }    \Bigg\}
+ a_1\frac{\exp\{-0.5 a_2 n^{ \alpha\, \tau }\} }{1-\Phi \left( x\right)}  \nonumber \\
&=& \exp\Bigg\{\theta_2 c_{2,\rho}  \frac{ (1+x)^{2+\rho }}{n^{\rho(\frac12 - \alpha)}  }    \Bigg\} \nonumber
\end{eqnarray}
uniformly for $0  \leq x   =o( \min\{ n^{\frac12 - \alpha}, \,  n^{\alpha   \tau /2} \}  ). $

\section{Proof of Theorem \ref{th3.3}  }\label{secend}
\setcounter{equation}{0}
We only give a proof for the case of $\rho \in (0, 1).$ The proof for the case of $\rho=1 $ is similar to the case of $\rho \in (0, 1).$
In the proof of theorem, we  use  the following lemma. The proof of the lemma is similar to the proof of Theorem A.6 of Hall and Heyde \cite{HH80}.
\begin{lemma}\label{HHA6}
Suppose that $X$ and $Y$ are random variables which are $\mathcal{F}_{j+n}^{\infty}$- and $\mathcal{F}_{j}$-measurable, respectively, and that $\mathbf{E}|X|^p < \infty,$ $\mathbf{E}|Y|^q  < \infty$,  where $p, q>1, p^{-1}+q^{-1}=1.$
Then
$$\Big|\mathbf{E}  X  Y  - \mathbf{E}  X\mathbf{E} Y \Big|\leq 2 [\psi(n)]^{1/p}   \big(\mathbf{E}|X|^p \big)^{1/p}  \big(\mathbf{E}|Y|^q \big)^{1/q} .$$
\end{lemma}

Denote by $\mathcal{F}_{l} = \sigma \{\eta_{i}, 1\leq i \leq  2ml-m  \}.   $ Then $Y_j$ is $\mathcal{F}_{j}$-measurable.
Since $\mathbf{E} \eta_{i } =0$ for all $i$, it is easy to see that for $1\leq j \leq k,$
\begin{eqnarray*}
 \Big|\mathbf{E}[ Y_j | \mathcal{F}_{ j-1 } ]\Big|&=& \Big|\sum_{i=1 }^{ m  } \Big( \mathbf{E}[\eta_{2m(j-1)+i } | \mathcal{F}_{j-1} ]- \mathbf{E} \eta_{2m(j-1)+i }  \Big)\Big| \\
  &\leq& \sum_{i=1 }^{ m  } \psi(m+i) \mathbf{E} |\eta_{2m(j-1)+i }|     \\
  &\leq& \sum_{i=1 }^{ m  } \psi(m+i) (\mathbf{E} |\eta_{2m(j-1)+i }|^{2+\rho})^{1/(2+\rho)}     \\
   &\leq& \sum_{i=1 }^{ m  } \psi(m+i) c_1,
\end{eqnarray*}
where $c_1$ is defined in (\ref{cosf3.11}).
Thus
\begin{eqnarray*}
\Big|\sum_{j=1}^k \mathbf{E}[ Y_j | \mathcal{F}_{ j-1 } ]  \Big| \leq  c_1\sum_{j=1}^k \sum_{i=1 }^{ m  } \psi(m+i) \leq n \psi(m ) c_1  .
\end{eqnarray*}
By (\ref{cosf3.11}),  we have
\begin{eqnarray}
 \mathbf{E}[  |Y_j-\mathbf{E}[ Y_j | \mathcal{F}_{ j-1 } ] |^{2+\rho} | \mathcal{F}_{ j-1 } ]&\leq& 2^{1+\rho}\mathbf{E}[  |Y_j|^{2+\rho} + |\mathbf{E}[ Y_j | \mathcal{F}_{ j-1 } ] |^{2+\rho} | \mathcal{F}_{ j-1 } ] \nonumber  \\
 &\leq& 2^{2+\rho}\mathbf{E}[  |Y_j|^{2+\rho}   | \mathcal{F}_{ j-1 } ] \nonumber \\
 &\leq& 2^{2+\rho} (1+ \psi(m )) \mathbf{E}  |Y_j|^{2+\rho} \nonumber  \\
  &\leq& 2^{2+\rho} (1+ \psi(m )) m^{1+\rho/2}c_1^{2+\rho}. \label{ine11.1}
\end{eqnarray}
Notice that  $\tau_n \rightarrow 0$ implies that $ m \psi^2(m ) \rightarrow 0$ as $n\rightarrow \infty$.
Similarly, by (\ref{cosdf3.12}), it holds
\begin{eqnarray}
  \mathbf{E}[  (Y_j-\mathbf{E}[ Y_j | \mathcal{F}_{ j-1 } ] )^2 | \mathcal{F}_{ j-1 } ]
  &=&   \mathbf{E}[  Y_j^2 | \mathcal{F}_{ j-1 } ]    -  ( \mathbf{E}[ Y_j | \mathcal{F}_{ j-1 } ] )^2  \nonumber \\
 &\geq& (1- \psi(m )) \mathbf{E}  Y_j^2   - ( \mathbf{E}[ Y_j | \mathcal{F}_{ j-1 } ] )^2 \nonumber \\
&\geq&   \frac12   c_2^2 m .\label{ine11.2}
\end{eqnarray}
Combining (\ref{ine11.1}) and (\ref{ine11.2}), we deduce that
$$\sum_{j=1}^k \mathbf{E}[  |Y_j-\mathbf{E}[ Y_j | \mathcal{F}_{ j-1 } ] |^2 | \mathcal{F}_{ j-1 } ] \asymp  n,$$
$$ \mathbf{E}[  |Y_j-\mathbf{E}[ Y_j | \mathcal{F}_{ j-1 } ] |^{2+\rho} | \mathcal{F}_{ j-1 } ]\leq c_\rho\, m^{\rho/2}  \mathbf{E}[  (Y_j-\mathbf{E}[ Y_j | \mathcal{F}_{ j-1 } ] )^2 | \mathcal{F}_{ j-1 } ]  $$
and, by Lemma \ref{HHA6},
\begin{eqnarray*}
&&\Big|\sum_{j=1}^k \mathbf{E}[  (Y_j-\mathbf{E}[ Y_j | \mathcal{F}_{ j-1 } ] )^2 | \mathcal{F}_{ j-1 } ] - \mathbf{E}  S_n^2 \Big| \\
 &&\leq\Big|\sum_{j=1}^k \mathbf{E}[  (Y_j-\mathbf{E}[ Y_j | \mathcal{F}_{ j-1 } ] )^2 | \mathcal{F}_{ j-1 } ] -\sum_{j=1}^k\mathbf{E}  Y_j^2 \Big| +\Big|  \mathbf{E}  S_n^2-\sum_{j=1}^k\mathbf{E}  Y_j^2 \Big|\\
  &&\leq  \sum_{j=1}^k \Big| \mathbf{E}[  Y_j^2 | \mathcal{F}_{ j-1 } ] -\mathbf{E}  Y_j^2 \Big|     + \sum_{j=1}^k \Big| \mathbf{E}[ Y_j | \mathcal{F}_{ j-1 } ] \Big|^2 +  \sum_{j\neq l}\Big|\mathbf{E}  Y_j  Y_l  \Big|\\
 &&\leq k\psi(m )\mathbf{E}Y_j^2 + k \Big | \sum_{i=1 }^{ m  } \psi(m+i) c_1 \Big |^2+ 2 \psi(m )^{1/2} \sum_{j\neq l} \sqrt{\mathbf{E}Y_j^2} \sqrt{\mathbf{E}Y_l^2}  \\
 &&\leq 2 n \psi(m )c_1^2+n m \psi^2(m ) c_1^2 + 2n \psi(m )^{1/2}k  c_1^2 .
\end{eqnarray*}
Denote by
$$\epsilon_n^2=  \psi(m ) + m \psi^2(m )   +   k\psi(m )^{1/2}  .$$
Applying  Theorem \ref{co0}   to $X_n:=\sum_{j=1}^k (Y_j-\mathbf{E}[ Y_j | \mathcal{F}_{ j-1 } ])/\sqrt{ \mathbf{E}  S_n^2}$,
we have for all $0\leq x =o(n^{\frac1 2-\alpha})  ,$
\begin{equation}
\bigg|\ln \frac{\mathbf{P}(X_n>x)}{1-\Phi \left( x\right)} \bigg| \leq  c_{ \rho} \,  \bigg(  \frac{ (1+x)^{2+\rho }}{n^{\rho(\frac12 - \alpha)}  }  + x^2 \epsilon_n^2 + (1+ x) \Big( \frac{1}{n^{\rho(\frac12 - \alpha)}  } +  \epsilon_n \Big)\bigg) .
\end{equation}
Notice that for $x\geq 0$ and $|\varepsilon|  \leq 1$,
\[
\frac{1-\Phi \left( x + \varepsilon \right)}{1-\Phi \left( x\right) }=   \exp\Big\{ O(1)(1+x) |\varepsilon| \Big\}
\]
and
\begin{eqnarray*}
\Big|\frac1{\sqrt{n} }\sum_{j=1}^k \mathbf{E}[ Y_j | \mathcal{F}_{ j-1 } ]  \Big| \leq     \sqrt{n} \psi(m ) c_1  .
\end{eqnarray*}
Thus
\begin{equation}
\bigg| \ln \frac{\mathbf{P}(S_n/\sqrt{ \mathbf{E}S_n^2}>x)}{1-\Phi \left( x\right)} \bigg| \leq  c_{\rho} \,  \bigg(  \frac{ (1+x)^{2+\rho }}{n^{\rho(\frac12 - \alpha)}  }  + x^2 \tau_n^2 + (1+ x) \Big( \frac{1}{n^{\rho(\frac12 - \alpha)}  } +  \tau_n  \Big)\bigg)  ,
\end{equation}
where $\tau_n^2$ is defined by (\ref{detan}).

\subsection*{Acknowledgements}
We would like to thank the associated editor for his helpful suggestions on the structure of paper.
We also deeply indebted to the anonymous referees  for their helpful comments.


\end{document}